\documentclass[11pt,reqno]{amsart}
\usepackage{amsmath,amsfonts, amssymb, amsthm}
  \usepackage{paralist}
  \usepackage{graphics} %% add this and next lines if pictures should be in esp format
  \usepackage{epsfig} %For pictures: screened artwork should be set up with an 85 or 100 line screen
\usepackage{graphicx}  \usepackage{epstopdf}%This is to transfer .eps figure to .pdf figure; please compile your paper using PDFLeTex or PDFTeXify.
 \usepackage[colorlinks=true]{hyperref}
\hypersetup{urlcolor=blue, citecolor=red}

\newtheorem{lemma}{Lemma}[section]
\newtheorem{theorem}{Theorem}[section]

\newtheorem{conjecture}{Conjecture}[section]

\theoremstyle{definition}

\newtheorem{remark}{Remark}[section]

\def\pmod #1{\ ({\rm{mod}}\ #1)}
\def\Z{\Bbb Z}
\def\N{\Bbb N}

\def\Q{\Bbb Q}

\def\l{\left}
\def\r{\right}
\def\bg{\bigg}
\def\({\bg(}
\def\){\bg)}
\def\t{\text}
\def\f{\frac}
\def\mo{{\rm{mod}\ }}
\def\pmod#1{\ (\mo\ #1)}
\def\Domb{\mathrm{Domb}}

\def\ls{\leq}

\def\sm{\setminus}
\def\bi{\binom}
\def\al{\alpha}

\def\eq{\equiv}

\def\Proof{\noindent{\it Proof}}

\begin{document}
%\hbox{Preprint, {\tt arXiv:2110.03651}}
\hbox{J. Comb. Number Theory, in press.}
\medskip

\title[New type series for powers of $\pi$]
%Use the shortened version of the full title
      {New type series for powers of $\pi$}
\author[Zhi-Wei Sun]{Zhi-Wei Sun}

\address{Department of Mathematics, Nanjing
University, Nanjing 210093, People's Republic of China}
\email{zwsun@nju.edu.cn}

\keywords{Ramanujan-type series, binomial coefficients, combinatorial identities.
\newline \indent 2020 {\it Mathematics Subject Classification}. Primary 11B65, 05A19; Secondary 33F10.
\newline \indent Supported by the Natural Science Foundation of China (grant no. 11971222).}

%The abstract of your paper
\begin{abstract} Motivated by Ramanujan-type series and Zeilberger-type series,
in this paper we investigate two new types of series for powers of $\pi$.
For example, we prove that
$$\sum_{k=0}^\infty(198k^2-425k+210)\frac{k^3\binom{2k}k^3}{4096^k}=-\frac1{21\pi}$$
and
$$\sum_{k=0}^\infty\frac{198k^2-227k+47}{\binom{2k}k^3}=\frac{3264-4\pi^2}{63}.$$
We also pose many conjectures in this new direction.
\end{abstract}
\maketitle

\section{Introduction}

In 1859, G. Bauer got the following identity:
$$\sum_{k=0}^\infty\f{4k+1}{(-64)^k}\bi{2k}k^3=\f2{\pi}.$$
This began the history of series for $\pi^{-1}$. In 1914, S. Ramanujan \cite{R}
posed 16 conjectural series of the following form:
\[\sum_{k=0}^\infty\f{bk+c}{m^k}a(k)=\f{\lambda\sqrt d}{\pi}\]
where $b,c,m\in\Z$ with $bm\not=0$, $d\in\Z^+=\{1,2,3,\ldots\}$ is squarefree,
 $\lambda\in\Q\sm\{0\}$, and $a(k)\in\Z$ is one of the products
\begin{equation}\label{term}\bi{2k}k^3,\ \bi{2k}k^2\bi{3k}k,\ \bi{2k}k^2\bi{4k}{2k},\ \bi{2k}k\bi{3k}k\bi{6k}{3k}.
\end{equation}
There are totally $36$ series of the above type, which are called
rational Ramanujan-type series for $\pi^{-1}$.
For surveys of such series, one may consult \cite{BB,Be,CC}
and S. Cooper \cite[Chapter 14]{Co17}). In general, a Ramanujan-type series (of type $R$)
has the form
 \[\sum_{k=0}^\infty\f{P_d(k)}{m^k}a(k)=\f{\al}{\pi^d},\tag{R}\]
where $P_d(x)\in\Z[x]$ with $\deg P=d>0$, $m$ is a nonzero integer, $\al$ is an algebraic number,
and $a(k)\in\Z$ for all $k\in\N=\{0,1,2,\ldots\}$. For example,
B. Gourevich used the PSLQ algorithm to find the conjectural identity
\begin{equation}\label{Go}\sum_{k=0}^\infty\f{168k^3+76k^2+14k+1}{2^{20k}}\bi{2k}k^7=\f{32}{\pi^3}.
\end{equation}

In this paper we introduce a new type of series different from Ramanujan's type.
A series of {\it type $S$} has the form
\[\sum_{k=0}^\infty\f{P(k)Q_{2d}(k)}{m^k}a(k)=\f{\al}{\pi^d},\tag{S}\]
where $P(x),Q_{2d}(x)\in\Z[x]$ with $\deg P>0$ and $\deg Q_{2d}\ls2d$, $m$ is a nonzero integer, $\al$
is an algebraic number, and
$a(k)\in\Z$ for all $k\in\N$.
We also investigate what happens if we replace  $\bi{2k}{k}$ in the general summand
of the classical Ramanujan-type series by $\bi{2k}{k+1}=kC_k$, where $C_k$ denotes the Catalan number
$\f1{k+1}\bi {2k}k=\bi{2k}k-\bi{2k}{k+1}$.

Now we state our first theorem which gives variants of the Ramanujan series involving cubes of central binomial coefficients.

\begin{theorem}\label{Th-A2} We have the following identities:
\begin{align*}\sum_{k=0}^\infty\f{(6k^2-19k+6)k^3\bi{2k}k^3}{256^k}&=-\f1{12\pi},\tag{S1}
\\\sum_{k=0}^\infty\f{(13608k^2+25050k+10589)\bi{2k}{k+1}^3}{256^k}&=27296-\f{84604}{\pi},\tag{S1$'$}
\\
\sum_{k=0}^\infty\f{(18k^2-29k+16)k^3\bi{2k}k^3}{(-512)^k}&=\f{\sqrt2}{24\pi},\tag{S2}
\\\sum_{k=0}^\infty\f{(83592k^2+152922k+64925)\bi{2k}{k+1}^3}{(-512)^k}&=
 625282\f{\sqrt2}{\pi}-281920,
\tag{S2$'$}
\\\sum_{k=0}^\infty\f{(198k^2-425k+210)k^3\bi{2k}k^3}{4096^k}&=-\f1{21\pi},\tag{S3}
\end{align*}
\[\begin{aligned}
&\sum_{k=0}^\infty\f{(32473224k^2+58012446k+24235261)\bi{2k}{k+1}^3}{4096^k}
\\&\qquad=667628032-\f{2097324016}{\pi}.
\end{aligned}
\tag{S3$'$}\]
\end{theorem}
\begin{remark} Motivated by \cite{vH}, we may also look at $p$-adic congruences corresponding to those identities in Theorem \ref{Th-A2}. For example, inspired by the identity (S1) we conjecture that for any prime $p>3$ we have
$$\sum_{k=1}^{(p-1)/2}\f{k^3\bi{2k}k^3}{256^k}(6k^2-19k+6)\eq(-1)^{\f{p+1}2}\f{p+2p^3}{48}\pmod{p^4}.$$
The paper \cite{S19} cotains many conjectural $p$-adic congruences.
\end{remark}

For our variants of other classical Ramanujan-type series, the reader may consult Section 4.

In 1993 D. Zeilberger \cite{Z} used the WZ method to establish the identity
\begin{equation}\label{Z-zeta(2)}\sum_{k=1}^\infty\f{21k-8}{k^3\bi{2k}k^3}=\f{\pi^2}6.
\end{equation}
More such series were conjectured by the author \cite{S11e}, and confirmed
by J. Guillera and M. Rogers \cite{GR}. In 1997, T. Amdeberhan and Zeilberger \cite{AZ}
used the WZ method to obtain the identity
\begin{equation}\label{AZ-zeta(3)}\sum_{k=1}^\infty\f{(-1)^k(205k^2-160k+32)}{k^5\bi{2k}k^5}=-2\zeta(3).
\end{equation}
In general, a Zeilberger-type series (of type $Z$) has the form
\[\sum_{k=1}^\infty\f{P_d(k)m^k}{k^{2d+1}a(k)}=cL(d+1,\chi)\tag{Z}\]
where $P_d(x)\in\Z[x]$ with $\deg P_d=d$, $m$ is a nonzero integer, $a(k)$ is a product of $2d+1$ binomial coefficients, $c$ is a nonzero rational number, and $L(2d,\chi)$ is the Dirichlet $L$ function associated with a Dirichlet character $\chi$.
For example, in 2010 the author conjectured that
\begin{equation}\label{Sun-zeta(3)}
\sum_{k=1}^\infty\f{(28k^2-18k+3)(-64)^k}{k^5\bi{2k}k^4\bi{3k}k}=-14\zeta(3)
\end{equation}
(cf. \cite{S13}) which remains open.

Motivated by the author's work \cite{S20} and the two new identities
$$
\sum_{k=1}^\infty\f{(22k^2-7k-3)64^k}{(2k+1)(3k+1)k^3\bi{2k}k^2\bi{3k}k}=48-4\pi^2$$
and
$$\sum_{k=1}^\infty\f{(22k^2+5k+1)64^k}{(2k+1)(3k+1)k^2\bi{2k}k^2\bi{3k}k}=4\pi^2-16
$$
missing from \cite[Theorem 1.2]{S20}, we introduce series of {\it type T}:
\[\sum_{k=1}^\infty\f{Q_{2d}(k)m^k}{a(k)}=c_0+c_1L(d+1,\chi),\tag{T}\]
where $Q_{2d}(x)\in\Z[x]$ with $\deg Q_{2d}\ls2d$, $m$ is a nonzero integer,
$a(k)$ is a product of $2d+1$ binomial coefficients, $c_0$ and $c_1$ are algebraic numbers,
and $\chi$ is a Dirichlet character. We also investigate what happens if we replace $\bi{2k}{k}$ in the denominator of a general summand
of the Zeilberger-type series (Z) by $\bi{2k}{k+1}=kC_k$.

Throught this paper, we adopt the notation
$$G=L\l(2,\l(\f{-4}{\cdot}\r)\r)=\sum_{k=0}^\infty\f{(-1)^k}{(2k+1)^2}$$
and $$K=L\l(2,\l(\f{-3}\cdot\r)\r)=\sum_{k=1}^\infty\f{(\f k3)}{k^2},$$
where $(\f d{\cdot})$ with $d\eq0,1\pmod4$ is the Kronecker symbol, and
$(\f k3)$ is the Legendre symbol. The number $G$ is usually called the Catalan constant.

Now we state our second theorem.

\begin{theorem} \label{Th-T}
We have the following identities:
\begin{align*}\sum_{k=0}^\infty\f{198k^2-227k+47}{\bi{2k}k^3}=&\f{3264-4\pi^2}{63},\tag{T1}
\\\sum_{k=1}^\infty\f{4884741k^2-7292783k+2041168}{\bi{2k}{k+1}^3}=&\f{68102148-10093535\pi^2}{126},\tag{T1$'$}
\\\sum_{k=0}^\infty\f{(54k^2-33k+18)(-8)^k}{\bi{2k}k^3}=&8+2G,\tag{T2}
\\\sum_{k=1}^\infty\f{(11277k^2-13124k+3212)(-8)^k}{\bi{2k}{k+1}^3}=&\f{38954G-24679}3,\tag{T2$'$}
\\\sum_{k=0}^\infty\f{(850k^2-1133k+69)8^k}{\bi{2k}k^2\bi{3k}k}=&\f{885-12\pi^2}5,\tag{T3}
\\\sum_{k=1}^\infty\f{(4621650k^2-7550827k+ 1785654)8^k}{\bi{2k}{k+1}^2\bi{3k}{k+1}}=&\f{6728880- 1468127\pi^2}5,\tag{T3$'$}
\\\sum_{k=0}^\infty\f{(18k^2-39k-6)16^k}{\bi{2k}k^3}=&\f{16-\pi^2}4,\tag{T4}
\\\sum_{k=1}^\infty\f{(2778k^2-6499k+1679)16^k}{\bi{2k}{k+1}^3}=&\f{89008-20305\pi^2}{12},\tag{T4$'$}
\end{align*}
\begin{align*}
\sum_{k=0}^\infty\f{(315k^2-54k+119)(-27)^k}{\bi{2k}k^2\bi{3k}k}=&\f{10+216K}5,\tag{T5}
\\\sum_{k=1}^\infty\f{(3391155k^2-3249747k+646792)(-27)^{k-1}}{\bi{2k}{k+1}^2\bi{3k}{k+1}}
=&\f{2344210-5115501K}{30},\tag{T5$'$}
\end{align*}
\begin{align*}
\sum_{k=0}^\infty\f{(1298k^2-7807k-3165)64^k}{\bi{2k}k^2\bi{3k}k}=&195-192\pi^2,\tag{T6}
\\\sum_{k=1}^\infty\f{(3470742k^2-20527615k+4631490)64^{k-1}}{\bi{2k}{k+1}^2\bi{3k}{k+1}}
=&387807- 198173\pi^2,\tag{T6$'$}
\\\sum_{k=0}^\infty\f{(2050k^2-5331k-1357)81^k}{\bi{2k}k^2\bi{4k}{2k}}=&\f{8800-2592\pi^2}{35},\tag{T7}
\end{align*}
\[\begin{aligned}&\sum_{k=1}^\infty\f{(17439700k^2 -47250409k+8929776)81^{k-1}}{\bi{2k}{k+1}^2\bi{4k}{2k+1}}
\\&\qquad=\f{43153685-44613764\pi^2}{210},\end{aligned}
\tag{T7$'$}\]
\begin{align*}
\sum_{k=0}^\infty\f{(50k^2+93k+22)(-144)^k}{\bi{2k}k^2\bi{4k}{2k}}=&\f{243K-200}4,\tag{T8}
\\\sum_{k=1}^\infty\f{(53725k^2-27494k+ 3675)(-144)^{k-1}}{\bi{2k}{k+1}^2\bi{4k}{2k+1}}=&\f{78296-122235K}{128}.\tag{T8$'$}
\end{align*}
\end{theorem}

Our third theorem is motivated by the Amdeberhan-Zeilberger identity \eqref{AZ-zeta(3)}.

\begin{theorem}\label{Th1.3} We have
\[\begin{aligned}&\sum_{k=0}^\infty\f{(-1)^k}{\bi{2k}k^5}
(684700k^4-1418358k^3+1100639k^2-365392k+47327)
\\&\qquad=\f{231168-64\zeta(3)}5.\end{aligned}\tag{T9}
\]
\end{theorem}

Recall that the Domb numbers are given by
$$\Domb(n)=\sum_{k=0}^n\bi nk^2\bi{2k}k\bi{2(n-k)}{n-k}\ \ (n=0,1,2,\ldots).$$
By H. H. Chan, S. H. Chan and Z. Liu \cite{CCL} and M. D. Rogers \cite{Rogers},
\begin{equation}\label{64-32}\sum_{k=0}^\infty\f{5k+1}{64^k}\Domb(k)=\f 8{\sqrt3\,\pi}
\ \t{and}\ \sum_{k=0}^\infty\f{3k+1}{(-32)^k}\Domb(k)=\f2{\pi}.\end{equation}
Franel numbers of order $4$ are those integers
$$f_n^{(4)}=\sum_{k=0}^n\bi nk^4\ \ (n\in\N).$$
In 2005 Yifan Yang discovered the first $\f1{\pi}$-series involving Franel numbers of order four:
\begin{equation}\label{Yang}\sum_{k=0}^\infty\f{4k+1}{36^k}f_k^{(4)}=\f{18}{\sqrt{15}\,\pi}.
\end{equation}
For more such series deduced via modular forms, see S. Cooper \cite{Co}.
Ap\'ery numbers of the second kind are those integers
$$\beta_n=\sum_{k=0}^n\bi nk^2\bi{n+k}k\ \ (n\in\N).$$
It is known (cf. Table 14.7 of \cite[p.\,653]{Co17}) that
\begin{gather*}\sum_{k=0}^\infty\f{5k+1}{72^k}\bi{2k}k\beta_k=\f{9\sqrt2}{2\pi},
\ \ \sum_{k=0}^\infty\f{11k+2}{147^k}\bi{2k}k\beta_k=\f{49\sqrt3}{10\pi},
\\ \sum_{k=0}^\infty\f{190k+29}{(-828)^k}\bi{2k}k\beta_k=\f{18\sqrt{23}}{\pi},
\ \ \sum_{k=0}^\infty\f{682k+71}{(-15228)^k}\bi{2k}k\beta_k=\f{162\sqrt{47}}{5\pi}.
\end{gather*}
For $n\in\N$ define the Cooper number
$$\mathrm{Co}(n)=\sum_{k=0}^n\bi nk^2\bi{n+k}k\bi{2k}n.$$
S. Cooper \cite{Co12} deduced the following three identities for $1/\pi$:
\begin{align}\sum_{k=0}^\infty\f{39k+10}{(-64)^k}\mathrm{Co}(k)&=\f{64}{\sqrt7\,\pi},
\\ \sum_{k=0}^\infty\f{21k+4}{125^k}\mathrm{Co}(k)&=\f{125}{8\pi},
\\\sum_{k=0}^\infty\f{11895k+1286}{(-22)^{3k}}\mathrm{Co}(k)&=\f{22^3}{\sqrt7\,\pi}.
\end{align}
Now we state our fourth theorem.

\begin{theorem}\label{Th1.4} {\rm (i)} We have
\begin{align}\label{D-32}\sum_{k=1}^\infty\f{k(3k+1)^2}{(-32)^k}\Domb(k)=-\f 2{\pi},
\\\label{D64}\sum_{k=1}^\infty\f{k^2(5k-3)}{64^k}\Domb(k)=\f{8\sqrt3}{9\pi}.
\end{align}

{\rm (ii)} We have the following identities:
\begin{align}\label{Y36}\sum_{k=1}^\infty\f{kf_k^{(4)}}{36^k}(40k^2-42k-1)=&\f{198\sqrt{15}}{25\pi},
\\\sum_{k=1}^\infty\f{kf_k^{(4)}}{(-64)^k}(80k^2+16k+3)=&-\f{544\sqrt{15}}{675\pi},
\\\sum_{k=1}^\infty\f{kf_k^{(4)}}{196^k}(33000k^2-3410k+291)=&\f{1582\sqrt{7}}{3\pi},
\\\sum_{k=1}^\infty\f{kf_k^{(4)}}{(-324)^k}(2720k^2+141k-7)=&-\f{3969\sqrt{5}}{200\pi},
\\\sum_{k=1}^\infty\f{kf_k^{(4)}}{1296^k}(156000k^2-2205k+563)=&\f{17901\sqrt{2}}{32\pi},
\\\sum_{k=1}^\infty\f{kf_k^{(4)}}{5776^k}(26079360k^2-81592k+33569)=&\f{220628\sqrt{95}}{75\pi}.
\end{align}

{\rm (iii)} We have \begin{align}\sum_{k=1}^\infty\f{k\bi{2k}k\beta_k}{72^k}(125k^2-300k-17)=&\f{999\sqrt{2}}{8\pi},
\\\sum_{k=1}^\infty\f{k\bi{2k}k\beta_k}{147^k}(121k^2-78k-1)=&\f{21903\sqrt{3}}{1250\pi},
\\\sum_{k=1}^\infty\f{k\bi{2k}k\beta_k}{(-828)^k}(5234500k^2+395850k-31627)=&\f{22302\sqrt{23}}{\pi},
\end{align}
and
\begin{align}\sum_{k=1}^\infty\f{k\bi{2k}k\beta_k}{(-15228)^k}(4127975500k^2-17838750k+6386881)
=-\f{3699918\sqrt{47}}{5\pi}.
\end{align}

{\rm (iv)} We have
\begin{align}\sum_{k=1}^\infty\f{k\mathrm{Co}(k)}{(-64)^k}(1365k^2+575k+86)&=-\f{7424\sqrt7}{147\pi},
\\\sum_{k=1}^\infty\f{k\mathrm{Co}(k)}{125^k}(1029k^2-413k+4)&=\f{16375}{96\pi},
\end{align}
and
\begin{equation}\sum_{k=1}^\infty\f{k\mathrm{Co}(k)}{(-22)^{3k}}(3480060675k^2+12712753k-5211590)=\f{226632032\sqrt7}
{147\pi}.
\end{equation}
\end{theorem}
\begin{remark} Actually our method to show Theorem \ref{Th1.4} allows us to deduce many more similar identities from all the known Ramanujan-type series for $1/\pi$ listed in \cite[pp.\,647--658]{Co17}.
\end{remark}

We will prove Theorem \ref{Th-A2} in Section 2, and Theorems \ref{Th-T}-\ref{Th1.3} in Section 3.
Section 4 contains our variants of classical Ramanujan-type series not included in Theorem \ref{Th-A2}.
Section 5 is devoted to our proof of Theorem \ref{Th1.4}.

Surprisingly, it seems that each non-classical Ramanujan-type series also has related series of type $S$. This enables us to pose many conjectural series of type $S$
in Sections 6-8 motivated by known series of type $R$.

Now we present our general conjecture relating series of type $R$ to series of type $S$.

\begin{conjecture} Suppose that we have a Ramanujan-type identity $(R)$
with $\limsup_{k\to+\infty}\root k\of{|a(k)|}<|m|$. For any non-constant polynomial $P(x)\in\Z[x]$, there is a polynomial $Q_{2d}(x)\in\Z$ of degree at most $2d$ such that
$$\sum_{k=0}^\infty\f{P(k)Q_{2d}(k)}{m^k}a(k)=\f{\al'}{\pi^d}$$
for some $\al'\in\Q(\al)$.
\end{conjecture}

Now we describe our algorithm to deduce series of type $S$ from a Ramanujan-type series (R)
with $a(k)$ a product of $2d+1$ binomial coefficients.
\medskip

\noindent {\bf General Algorithm}. Suppose that we have a Ramanujan-type series $(R)$
with $a(k)$ a product of $2d+1$ binomial coefficients.
Given a monic polynomial $P(x)\in\Z[x]$ with $\deg P>0$, we first write
$P(x)=(x-\al_1)\cdots(x-\al_r)$ with $\al_1,\ldots,\al_r$ complex numbers.
Using the Gosper algorithm (cf. \cite{PWZ}) we find a polynomial $Q_{2d+1}$
of degree at most $2d+1$ such that
$\sum_{k=0}^n\f{a(k)}{m^k}Q_{2d+1}(k)$ has a closed form which tends to $0$,
so that we get
\begin{equation}\label{Q}\sum_{k=0}^\infty\f{a(k)}{m^k}Q_{2d+1}(k)=0.
\end{equation}
Note that \eqref{Q}$\times P_d(\al_1)-(R)\times Q_{2d+1}(\al_1)$ yields
an identity of the form
\begin{equation}\label{R2d}\sum_{k=0}^\infty\f{(k-\al_1)a(k)}{m^k}R_{2d}(k)=-Q_{2d+1}(\al_1)\f{\al}{\pi^d},
\end{equation}
where $R_{2d}$ is a polynomial of degree at most $2d$. Using the Gosper algorithm, we find
\begin{equation}\label{Q1}\sum_{k=0}^\infty\f{(k-\al_1)a(k)}{m^k}\tilde Q_{2d+1}(k)=0,
\end{equation}
where $\tilde Q_{2d+1}$ is a polynomial of degree at most $2d+1$.
Then \eqref{Q1}$\times R_{2d}(\al_2)-\tilde Q_{2d+1}(\al_2)\times$\eqref{R2d} yields
an identity of the form
\begin{equation*}\label{R2d-}\sum_{k=0}^\infty\f{(k-\al_1)(k-\al_2)a(k)}{m^k}\tilde R_{2d}(k)=Q_{2d+1}(\al_1)\tilde Q_{2d+1}(\al_2)\f{\al}{\pi^d}
\end{equation*}
with $\tilde R_{2d}$ a polynomial of degree at most $2d$. Continue this process, we finally
obtain the exact value of
$$\sum_{k=0}^\infty\f{P(k)a(k)}{m^k}Q_{2d}(k)$$
with $Q_{2d}$ a suitable polynomial of degree at most $2d$.
\medskip

\begin{remark} In practice this algorithm works well. Actually it could also be used to deduce series of type $T$ from series of type $Z$.
\end{remark}
\medskip
\noindent
{\bf Example 1.1}. Let's implement the algorithm to deduce a series of type $S$ with $P(k)=k^2+1=(k+i)(k-i)$
from Ramanujan's identity
\begin{equation}\label{256}\sum_{k=0}^\infty\f{\bi{2k}k^3}{256^k}(6k+1)=\f4{\pi}.
\end{equation}
By Gosper's algorithm we find that
$$\sum_{k=0}^n\f{\bi{2k}k^3}{256^k}(24k^3-12k^2-6k-1)=-(2n+1)^2\f{\bi{2n}n^3}{256^n}.$$
Letting $n\to+\infty$ we get
\begin{equation}\label{dd}\sum_{k=0}^\infty\f{\bi{2k}k^3}{256^k}(24k^3-12k^2-6k-1)=0.
\end{equation}
Note that $6k+1=6(k+i)+1-6i$ and
$$24k^3-12k^2-6k-1=6(k+i)(4k^2-(2+4i)k+2i-5)+30i+11.$$
Via  \eqref{256}$\times (30i+11)+(6i-1)\times$\eqref{dd} we obtain
\begin{equation}\label{ee}
\sum_{k=0}^\infty\f{(k+i)\bi{2k}k^3}{256^k}((k-i)(6i-1)(2k-1)+8)=\f{11+30i}{3\pi}.
\end{equation}
Via the Gosper algorithm, we find that
\begin{equation}\label{ff}\sum_{k=0}^\infty\f{(k+i)\bi{2k}k^3}{256^k}
(2(k-i)((132+360i)k^2-(602+144i)k+15-708i)+1373)=0.
\end{equation}
Observe that \eqref{ee}$\times1373-8\times$\eqref{ff} yields
\begin{equation}\label{k^2+1-256}\sum_{k=0}^\infty\f{(k^2+1)\bi{2k}k^3}{256^k}(192k^2-626k-103)=-\f{1373}{3\pi}.
\end{equation}
We can also deduce other identities similar to \eqref{k^2+1-256} such as
\begin{align}\sum_{k=0}^\infty\f{(k^2+1)\bi{2k}k^3}{(-512)^k}(6k+1)=&\f{11\sqrt2}{6\pi},
\\\sum_{k=0}^\infty\f{(k^2+1)\bi{2k}k^3}{4096^k}(126504k^2-921334k-109205)=-&\f{1063412}{3\pi}.
\end{align}

\section{Proof of Theorem \ref{Th-A2}}

\setcounter{theorem}{0}
\setcounter{equation}{0}
\setcounter{conjecture}{0}

\begin{lemma}\label{Lem-222} Let $m$ be any nonzero integer, and let $n\in\N=\{0,1,2,\ldots\}$.

{\rm (i)} We have
\begin{gather*}\sum_{k=0}^n\f{\bi{2k}k^3}{m^k}\l((64-m)k^3+96k^2+48k+8\r)
=8(2n+1)^3\f{\bi{2n}n^3}{m^n},
\\\sum_{k=0}^n\f{k\bi{2k}k^3}{m^k}\l((64-m)k^3+(96+m)k^2+48k+8\r)
=8n(2n+1)^3\f{\bi{2n}n^3}{m^n},
\end{gather*}
and
\begin{equation*}\begin{aligned}
&\sum_{k=0}^n\f{k^2\bi{2k}k^3}{m^k}\l((64-m)k^3+(96+2m)k^2+(48-m)k+8\r)
\\&\qquad\qquad=8n^2(2n+1)^3\f{\bi{2n}n^3}{m^n}.
\end{aligned}\end{equation*}

{\rm (ii)} We have
\begin{align*}&\sum_{k=0}^n\f{\bi{2k}k^3}{(k+1)^3m^k}\l((64-m)k^3+(96-3m)k^2+(48-3m)k+8-m\r)
\\=&-m+8(2n+1)^3\f{\bi{2n}n^3}{(n+1)^3m^n},
\end{align*}
and
\begin{align*}
&\sum_{k=0}^n\f{k\bi{2k}k^3}{(k+1)^3m^k}P(k,m)
=8m-8(2n+1)^3(m-8n+mn)\f{\bi{2n}n^3}{(n+1)^3m^n},
\end{align*}
where $P(k,m)$ denotes
$$(512-72m+m^2)k^3+(768-176m+3m^2)k^2
+(384-144m+3m^2)k+64-40m+m^2.$$
Also,
\begin{align*}
&\sum_{k=0}^n\f{k^2\bi{2k}k^3}{(k+1)^3m^k}Q(k,m)+8m(m+8)
\\=&8(2n+1)^3(m(m+8)-40mn(n+1)+2m^2n+64n^2+m^2n^2)
\f{\bi{2n}n^3}{(n+1)^3m^n},
\end{align*}
where $Q(k,m)$ denotes
\begin{align*}&(4096-2624m+1024m^2-m^3)k^3+(6144-6464m+304m^2-3m^3)k^2
\\&+(3072-5120m+296m^2-3m^3)k+512-1344m+96m^2-m^3.
\end{align*}

\end{lemma}
\Proof. It is easy to prove all the six identities by induction on $n$. We find them via the Gosper algorithm (cf. \cite{PWZ}). \qed

\begin{lemma} [Sun \cite{S20}] We have
\begin{align}\sum_{k=0}^\infty\f{k(6k-1)\bi{2k}k^3}{(2k-1)^3256^k}&=\f1{2\pi},
\\\sum_{k=0}^\infty\f{(30k^2+3k-2)\bi{2k}k^3}{(2k-1)^3(-512)^k}&=\f{27\sqrt2}{8\pi},
\\\label{4096(2k-1)}\sum_{k=0}^\infty\f{(42k^2-3k-1)\bi{2k}k^3}{(2k-1)^3 4096^k}&=\f{27}{8\pi}.
\end{align}
\end{lemma}

\medskip
\noindent {\it Proof of Theorem \ref{Th-A2}}. By Stirling's formula,
$$\bi{2n}n\sim\f{4^n}{\sqrt{n\pi}}\ \ \t{as}\ n\to+\infty.$$
Thus, for any integer $m$ with $|m|>64$, by Lemma \ref{Lem-222} we have
\begin{align}\label{k^0}\sum_{k=0}^\infty\f{\bi{2k}k^3}{m^k}\l((64-m)k^3+96k^2+48k+8\r)
=&0,
\\\label{k}\sum_{k=0}^\infty\f{k\bi{2k}k^3}{m^k}\l((64-m)k^3+(96+m)k^2+48k+8\r)
=&0,
\\\label{k^2}\sum_{k=0}^n\f{k^2\bi{2k}k^3}{m^k}\l((64-m)k^3+(96+2m)k^2+(48-m)k+8\r)=&0;
\end{align}
also
\begin{equation}\label{(k+1)^0}\sum_{k=0}^n\f{\bi{2k}k^3}{(k+1)^3m^k}\l((64-m)k^3+(96-3m)k^2+(48-3m)k+8\r)
=-m,\end{equation}
\begin{equation}\label{(k+1)}\sum_{k=0}^n\f{k\bi{2k}k^3}{(k+1)^3m^k}P(k,m)
=8m\end{equation}
and
\begin{equation}\label{(k+1)^2}\sum_{k=0}^n\f{k^2\bi{2k}k^3}{(k+1)^3m^k}Q(k,m),
=-8m(m+8)\end{equation}
where $P(k,m)$ and $Q(k,m)$ are as in Lemma \ref{Lem-222}(ii).

Now we prove (S3) in details. (The identities (S1) and (S2) can be proved similarly.)
Taking $m=4096$ in \eqref{k^0}-\eqref{k^2}, we get the identities
\begin{align}\label{4096-0}\sum_{k=0}^\infty\f{\bi{2k}k^3}{4096^k}\l(504k^3-12k^2-6k-1\r)=&0,
\\\label{4096-1}\sum_{k=0}^\infty\f{k\bi{2k}k^3}{4096^k}\l(504k^3-524k^2-6k-1\r)=&0,
\\\label{4096-2}\sum_{k=0}^\infty\f{k^2\bi{2k}k^3}{4096^k}\l(504k^3-1036k^2+506k-1\r)=&0.
\end{align}
Combining the Ramanujan identity
$$\sum_{k=0}^\infty\f{\bi{2k}k^3}{4096^k}(42k+5)=\f{16}{\pi}$$
with \eqref{4096-0}, and noting that
$$12k(210k^2-5k+1)=60k(42k^2-k)+12k = 5\times 12k(42k^2-k)+12k,$$
we obtain
$$\sum_{k=0}^\infty\f{\bi{2k}k^3}{4096^k}12k(210k^2-5k+1)=
\sum_{k=0}^\infty\f{\bi{2k}k^3}{4096^k}(5(6k+1)+12k)=\f{16}{\pi}$$
and hence
$$\sum_{k=0}^\infty\f{k\bi{2k}k^3}{4096^k}(210k^2-5k+1)=\f4{3\pi}.$$
Adding this and \eqref{4096-1}, we see that
$$\sum_{k=0}^\infty\f{k^2\bi{2k}k^3}{4096^k}(504k^2-314k-11)=\f4{3\pi}.$$
\eqref{4096-2} times $11$ minus the last formula yields
$$\sum_{k=0}^\infty\f{k^2\bi{2k}k^3}{4096^k}(11\times504k^3-(11\times1036+504)k^2+(11\times506+314)k)=-\f4{3\pi},$$
which is equivalent to (S3).

Next we prove (S3$'$) in details. (The identities (S1$'$) and (S2$'$) can be proved similarly.)
For any positive integer $k$, clearly
$$\f{\bi{2k}k}{2k-1}=\f2k\bi{2(k-1)}{k-1}.$$
Thus, by \eqref{4096(2k-1)} we have
\begin{align*}\f{27}{8\pi}-1=&\sum_{k=1}^\infty\f{42k^2-3k-1}{4096^k}\l(\f 2k\bi{2(k-1)}{k-1}\r)^3
\\=&\f 8{4096}\sum_{k=0}^\infty\f{42(k+1)^2-3(k+1)-1}{4096^k}\cdot\f{\bi{2k}k^3}{(k+1)^3}
\end{align*}
and hence
\begin{equation}\label{(k+1)^0-4096}\sum_{k=0}^\infty\f{\bi{2k}k^3}{(k+1)^34096^k}(42k^2+81k+38)=\f{1728}{\pi}-512.
\end{equation}
On the other hand, applying \eqref{(k+1)^0} with $m=4096$, we get
\begin{equation}\label{(k+1)^0-4096'}\sum_{k=0}^\infty\f{\bi{2k}k^3}{(k+1)^34096^k}
(504k^3+1524k^2+1530k+511)=512.
\end{equation}
Note that \eqref{(k+1)^0-4096'}$\times 38-$\eqref{(k+1)^0-4096}$\times511$ yields
\begin{equation}\label{(k+1)-4096}
\sum_{k=0}^\infty\f{k\bi{2k}k^3}{(k+1)^34096^k}(2128k^2+4050k+1861)=31232-\f{98112}{\pi}.
\end{equation}
Taking $m=4096$ in \eqref{(k+1)}, we obtain
\begin{equation}\label{(k+1)-4096'}\sum_{k=0}^\infty\f{k\bi{2k}k^3}{(k+1)^34096^k}
(257544k^3+775180k^2+777222k+259585)=512.
\end{equation}
Observe that \eqref{(k+1)-4096'}$\times 1861-$\eqref{(k+1)-4096}$\times259585$ gives
\begin{equation}\label{(k+1)^2-4096}
\sum_{k=0}^\infty\f{k^2\bi{2k}k^3}{(k+1)^34096^k}(78162k^2+145175k+64431)
=\f{4153360}{\pi}-1321984.
\end{equation}
Putting $m=4096$ in \eqref{(k+1)^2}, we get
\begin{equation}\label{(k+1)^2-4096'}\begin{aligned}&\sum_{k=0}^\infty\f{k^2\bi{2k}k^3}{(k+1)^34096^k}
(130830840k^3+392743412k^2+392994810k+131082751)
\\&\qquad\qquad=262656.\end{aligned}
\end{equation}
Note that \eqref{(k+1)^2-4096'}$\times 64431-$\eqref{(k+1)^2-4096}$\times131082751$ yields
\begin{align*}&\sum_{k=0}^\infty\f{k^3\bi{2k}k^3}{(k+1)^34096^k}
(32473224k^2+58012446k+24235261)
\\&\quad=667628032-\f{2097324016}{\pi},
\end{align*}
which is equivalent to the desired (S3') since $\bi{2k}{k+1}=\f k{k+1}\bi{2k}k$.
 \qed

 \section{Proofs of Theorems \ref{Th-T} and \ref{Th1.3}}

 \setcounter{theorem}{0}
\setcounter{equation}{0}
\setcounter{conjecture}{0}

\begin{lemma} \label{Lem-Z} Let $m\in\Z$ and $n\in\Z^+=\{1,2,3,\ldots\}$.

{\rm (i)} We have the following identities:
\begin{align*} \sum_{k=1}^n\f{m^k}{k\bi{2k}k^3}
((m-64)k^3+(2m+96)k^2-48k+8)=&-m+\f{(n+1)^2m^{n+1}}{\bi{2n}n^3},
\\\sum_{k=1}^n\f{m^k}{k^2\bi{2k}k^3}
((m-64)k^3+(m+96)k^2-48k+8)=&-m+\f{(n+1)m^{n+1}}{\bi{2n}n^3},
\\\sum_{k=1}^n\f{m^k}{k^3\bi{2k}k^3}
((m-64)k^3+96k^2-48k+8)=&-m+\f{m^{n+1}}{\bi{2n}n^3},
\end{align*}

{\rm (ii)} We have
\begin{align*}\sum_{k=1}^n\f{(k+1)m^k}{k^3\bi{2k}k^3}
P_1(k,m)=m(m-432)+(216(n+2)-m(n+1))\f{m^{n+1}}{\bi{2n}n^3},
\end{align*}
\begin{align*}&\sum_{k=1}^n\f{(k+1)^2m^k}{k^3\bi{2k}k^3}
P_2(k,m)
-m(m^2-216m+186624)
\\=&-\f{m^{n+1}}{\bi{2n}n^3}
\\&\times(186624(n+1)-216m(n^2+1)+m^2(2n+1)-648mn+m^2n^2+46656n^2),
\end{align*}
where
\begin{align*}P_1(k,m)=&(280m-m^2-13824)k^3+(56m+20736)k^2
\\&-(8m+10368)k+1728
\end{align*}
and
\begin{align*}P_2(k,m)=&(2985984-60480m+280m^2-m^3)k^3
\\&-(4478976+58752m+8m^2)k^2+(2239488+15552m)k
\\&-(1728m+373248).
\end{align*}
\end{lemma}
\Proof. It is easy to prove the five identities by induction on $n$. We find them via the Gosper algorithm (cf. \cite{PWZ}). \qed

\medskip
\noindent {\it Proof of {\rm (T1)}}. Lemma \ref{Lem-Z}(i) with $m=1$ gives the identities
\begin{align*}\sum_{k=1}^n\f{63k^3-98k^2+47k-8}{k\bi{2k}k^3}=&1-\f{(n+1)^2}{\bi{2n}n^3},
\\\sum_{k=1}^n\f{63k^3-97k^2+48k-8}{k^2\bi{2k}k^3}=&1-\f{n+1}{\bi{2n}n^3},
\\\sum_{k=1}^n\f{63k^3-96k^2+48k-8}{k^3\bi{2k}k^3}=&1-\f{1}{\bi{2n}n^3}.
\end{align*}
Letting $n\to+\infty$ we get
\begin{align}\label{1-1}\sum_{k=1}^\infty\f{63k^3-98k^2+47k-8}{k\bi{2k}k^3}=&1,
\\\label{1-2}\sum_{k=1}^\infty\f{63k^3-97k^2+48k-8}{k^2\bi{2k}k^3}=&1,
\\\label{1-3}\sum_{k=1}^\infty\f{63k^3-96k^2+48k-8}{k^3\bi{2k}k^3}=&1.
\end{align}
Note that \eqref{1-3} minus \eqref{Z-zeta(2)} yields the identity
\begin{equation}\label{k=2}\sum_{k=1}^\infty\f{21k^2-32k+9}{k^2\bi{2k}k^3}=\f13-\f{\pi^2}{18}.
\end{equation}
Via \eqref{1-2}$\times9+8\times$\eqref{k=2}, we obtain
\begin{equation}\label{k=1}\sum_{k=1}^\infty\f{567k^2-705k+176}{k\bi{2k}k}=\f{35}3-\f49\pi^2.
\end{equation}
Observe that \eqref{1-1}$\times22+$\eqref{k=1} yields the identity
$$\sum_{k=1}^\infty\f{198k^2-227k+47}{\bi{2k}k^3}=\f{303-4\pi^2}{63},$$
which is equivalent to (T1) since $303+47\times63=3264$. \qed

\medskip
\noindent{\it Proof of {\rm (T1$'$)}}.  Lemma \ref{Lem-Z}(ii) with $m=1$ gives the identities
\begin{align*}&\sum_{k=1}^n\f{k+1}{k^3\bi{2k}k^3}(13545k^3-20792k^2+10376k-1728)
\\&\qquad=431-\f{215n+431}{\bi{2n}n^3}
\end{align*}
and
\begin{align*}
&\sum_{k=1}^n\f{(k+1)^2}{k^3\bi{2k}k^3}(2925783k^2-4537736k^2+2255040k-374976)
\\&\qquad=186409-\f{46441n^2+185978n+186409}{\bi{2n}n^3}.
\end{align*}
Letting $n\to+\infty$ we get
\begin{align*}\sum_{k=1}^n\f{k+1}{k^3\bi{2k}k^3}(13545k^3-20792k^2+10376k-1728)=&431,
\\\sum_{k=1}^n\f{(k+1)^2}{k^3\bi{2k}k^3}(2925783k^3-4537736k^2+2255040k-374976)
=&186409.
\end{align*} They can be rewritten as follows:
\begin{align}\label{r-1}\sum_{k=1}^n\f{k+1}{k^3\bi{2k}k^3}((k+1)(13545k^2-34337k+44713)-46441)=&431,
\\\label{r-2}\sum_{k=1}^n\f{(k+1)^2}{k^3\bi{2k}k^3}((k+1)(2925783k^2-7463519k+9718559)-10093535)
=&186409.
\end{align}

We may rewrite \eqref{1-3} in the form
\begin{equation}\label{r-0}
\sum_{k=1}^\infty\f{(k+1)(63k^2-159k+207)-215}{k^3\bi{2k}k^3}=1,
\end{equation}
and rewrite \eqref{Z-zeta(2)} in the form
\begin{equation}\label{r-0'}\sum_{k=1}^\infty\f{21(k+1)-29}{k^3\bi{2k}k^3}=\f{\pi^2}6.
\end{equation}
Note that \eqref{r-0}$\times29-215\times$\eqref{r-0'} yields
$$\sum_{k=1}^\infty\f{k+1}{k^3\bi{2k}k^3}(609k^2-1537k+496)=\f{29}3-\f{215}{18}\pi^2,$$
which has the equivalent form
\begin{equation}\label{r-1'}
\sum_{k=1}^\infty\f{k+1}{k^3\bi{2k}k^3}((k+1)(609k-2146)+2642)=\f{29}3-\f{215}{18}\pi^2.
\end{equation}
Via \eqref{r-1}$\times2642+46441\times$\eqref{r-1'} we obtain
\begin{align*}&\sum_{k=1}^\infty\f{(k+1)^2}{k^3\bi{2k}k^3}(166446 k^2- 290399 k+85904) =\f{22153}3-\f{46441}{18}\pi^2,
\end{align*} which has the equivalent form
\begin{equation}\label{r-2'}\begin{aligned}&\sum_{k=1}^\infty\f{(k+1)^2}{k^3\bi{2k}k^3}
((k+1)(166446 k- 456845)+542749)
 =\f{22153}3-\f{46441}{18}\pi^2.
\end{aligned}\end{equation}
Observe that \eqref{r-2}$\times542749+10093535\times$\eqref{r-2'} yields
$$\sum_{k=1}^\infty\f{(k+1)^3}{k^3\bi{2k}k^3}(2041168 - 7292783 k + 4884741 k^2)
=\f{68102148 - 10093535\pi^2}{126},$$
which is equivalent to the desired (T1$'$) since $\bi{2k}{k+1}=\f k{k+1}\bi{2k}k$.
This concludes the proof. \qed

\begin{remark}
In the spirit of our proofs of (T1) and (T1$'$), we can prove (T2),(T2$'$),(T4),(T4$'$) similarly
by using Lemma \ref{Lem-Z} and the known identities
$$\sum_{k=1}^\infty\f{(3k-1)(-8)^k}{k^3\bi{2k}k^3}=-2G
\ \t{and}\ \ \sum_{k=1}^\infty\f{(3k-1)16^k}{k^3\bi{2k}k^3}=\pi^2$$
(cf. \cite{GR}).
\end{remark}

\begin{lemma} \label{Lem-Z23} Let $m\in\Z$ and $n\in\Z^+$.

{\rm (i)} We have the following identities:
\begin{align*} &\sum_{k=1}^n\f{m^k}{k\bi{2k}k^2\bi{3k}k}
((m-108)k^3+(2m+162)k^2+(m-78)k+12)
\\&\qquad=-m+\f{(n+1)^2m^{n+1}}{\bi{2n}n^2\bi{3n}n},
\\&\sum_{k=1}^n\f{m^k}{k^2\bi{2k}k^2\bi{3k}k}
((m-108)k^3+(m+162)k^2-78k+12)
\\&\qquad=-m+\f{(n+1)m^{n+1}}{\bi{2n}n^2\bi{3n}n},
\\&\sum_{k=1}^n\f{m^k}{k^3\bi{2k}k^2\bi{3k}k}
((m-108)k^3+162k^2-78k+12)
\\&\qquad=-m+\f{m^{n+1}}{\bi{2n}n^2\bi{3n}n}.
\end{align*}

{\rm (ii)} We have
\begin{align*}\sum_{k=1}^n\f{(k+1)m^k}{k^3\bi{2k}k^2\bi{3k}k}
P_1(k,m)=m(m-720)+(360(n+2)-m(n+1))\f{m^{n+1}}{\bi{2n}n^2\bi{3n}n}
\end{align*}
and
\begin{align*}&\sum_{k=1}^n\f{(k+1)^2m^k}{k^3\bi{2k}k^2\bi{3k}k}
P_2(k,m)
-m(m^2-372m+518400)
\\=&-\f{m^{n+1}}{\bi{2n}n^2\bi{3n}n}
\\&\times(m^2(n+1)^2-366mn(n+3)+129600n^2+518400(n+1)-372m),
\end{align*}
where
\begin{align*}P_1(k,m)=&(468m-m^2-38880)k^3+(90m+58320)k^2
\\&-(12m+28080)k+4320
\end{align*}
and
\begin{align*}P_2(k,m)=&(13996800-169128m+474m^2-m^3)k^3
\\&-(20995200+160380m+12m^2)k^2
\\&+(10108800+41112m)k-1555200
\end{align*}
\end{lemma}

\begin{lemma} \label{Lem-Z24} Let $m\in\Z$ and $n\in\Z^+$.

{\rm (i)} We have the following identities:
\begin{align*} &\sum_{k=1}^n\f{m^k}{k\bi{2k}k^2\bi{4k}{2k}}
((m-256)k^3+(2m+384)k^2+(m-176)k+24)
\\&\qquad=-m+\f{(n+1)^2m^{n+1}}{\bi{2n}n^2\bi{4n}{2n}},
\\&\sum_{k=1}^n\f{m^k}{k^2\bi{2k}k^2\bi{4k}{2k}}
((m-256)k^3+(m+384)k^2-176k+24)
\\&\qquad=-m+\f{(n+1)m^{n+1}}{\bi{2n}n^2\bi{4n}{2n}},
\\&\sum_{k=1}^n\f{m^k}{k^3\bi{2k}k^2\bi{4k}{2k}}
((m-256)k^3+384k^2-176k+24)
\\&\qquad=-m+\f{m^{n+1}}{\bi{2n}n^2\bi{4n}{2n}}.
\end{align*}

{\rm (ii)} We have
\begin{align*}\sum_{k=1}^n\f{(k+1)m^k}{k^3\bi{2k}k^2\bi{4k}{2k}}
P_1(k,m)=m(m-1680)+(840(n+2)-m(n+1))\f{m^{n+1}}{\bi{2n}n^2\bi{4n}{2n}}
\end{align*}
and
\begin{align*}&\sum_{k=1}^n\f{(k+1)^2m^k}{k^3\bi{2k}k^2\bi{4k}{2k}}
P_2(k,m)
-m(m^2-904m+2822400)
\\=&-\f{m^{n+1}}{\bi{2n}n^2\bi{4n}{2n}}
\\&\times(m^2(n+1)^2-872mn(n+3)+705600n(n+4)-904m+2822400),
\end{align*}
where
\begin{align*}P_1(k,m)=&(1096m-m^2-215040)k^3+(200m+322560)k^2
\\&-(24m+147840)k+20160
\end{align*}
and
\begin{align*}P_2(k,m)=&(180633600-928832m+1128m^2-m^3)k^3
\\&-(270950400+853120m+24m^2)k^2
\\&+(124185600+209088m)k-(16934400+20160m).
\end{align*}
\end{lemma}

Using Lemmas \ref{Lem-Z23} and \ref{Lem-Z24}
as well as some results in \cite{S20}, we can prove all the remaining formulas in Theorem \ref{Th-T}.

\medskip
\noindent{\it Proof of Theorem 1.3}. By Gosper's algorithm, we find the identity
$$\sum_{k=1}^n\f{(-1)^k}{k^5\bi{2k}k^5}(1025k^5-2560k^4+2560k^3-1280k^2+320k-32)=-1+\f{(-1)^n}{\bi{2n}n}$$
which can be easily proved by induction on $n$. Letting $n\to+\infty$ we get
\begin{equation*}\sum_{k=1}^\infty\f{(-1)^k}{k^5\bi{2k}k^5}(1025k^5-2560k^4+2560k^3-1280k^2+320k-32)=-1.
\end{equation*}
Adding this to \eqref{AZ-zeta(3)} we get
\begin{equation}
\label{k-4}\sum_{k=1}^\infty\f{(-1)^k}{k^4\bi{2k}k^5}(205k^4-512k^3+512k^2-215k+32)=-\f{1+2\zeta(3)}5.
\end{equation}

By Gosper's algorithm, we find that
$$\sum_{k=1}^n\f{(-1)^k}{k^4\bi{2k}k^5}(1025k^5-2559k^4+2560k^3-1280k^2+320k-32)=-1+\f{(-1)^n(n+1)}{\bi{2n}n}$$
and hence
\begin{equation}\label{k-4'}\sum_{k=1}^\infty\f{(-1)^k}{k^4\bi{2k}k^5}(1025k^5-2559k^4+2560k^3-1280k^2+320k-32)=-1.
\end{equation}
Adding \eqref{k-4} and \eqref{k-4'}, we obtain
\begin{equation}\label{k-3}\sum_{k=1}^\infty\f{(-1)^k}{k^3\bi{2k}k^5}
(1025k^4-2354k^3+2048k^2-768k+105)=-\f25(3+\zeta(3)).
\end{equation}

By Gosper's algorithm, we find that
$$\sum_{k=1}^n\f{(-1)^k}{k^3\bi{2k}k^5}(1025k^5-2558k^4+2561k^3-1280k^2+320k-32)=-1+\f{(-1)^n(n+1)^2}{\bi{2n}n}$$
and hence
\begin{equation}\label{k-3'}\sum_{k=1}^\infty\f{(-1)^k}{k^3\bi{2k}k^5}(1025k^5-2558k^4+2561k^3-1280k^2+320k-32)=-1.
\end{equation}
Note that \eqref{k-3}$\times32+105\times$\eqref{k-3'} yields
\begin{equation}\label{k-2}\begin{aligned}&\sum_{k=1}^\infty\f{(-1)^k}{k^2\bi{2k}k^5}
(107625k^4-235790k^3+193577k^2-68864k+9024)
\\&\qquad=-\f{64}5(3+\zeta(3))-105.\end{aligned}\end{equation}

By Gosper's algorithm, we find that
$$\sum_{k=1}^n\f{(-1)^k}{k^2\bi{2k}k^5}(1025k^5-2557k^4+2563k^3-1279k^2+320k-32)=-1+\f{(-1)^n(n+1)^3}{\bi{2n}n}$$
and hence
\begin{equation}\label{k-2'}\sum_{k=1}^\infty\f{(-1)^k}{k^2\bi{2k}k^5}(1025k^5-2557k^4+2563k^3-1279k^2+320k-32)=-1.
\end{equation}
Via \eqref{k-2}$+282\times$\eqref{k-2'} we get
\begin{equation}\label{k-1}\begin{aligned}&\sum_{k=1}^\infty\f{(-1)^k}{k\bi{2k}k^5}
(289050k^4-613449k^3+486976k^2-167101k+21376)
\\&\qquad=-\f{64}5(3+\zeta(3))-105-282.\end{aligned}\end{equation}

By Gosper's algorithm, we find that
$$\sum_{k=1}^n\f{(-1)^k}{k\bi{2k}k^5}(1025k^5-2556k^4+2566k^3-1276k^2+321k-32)=-1+\f{(-1)^n(n+1)^4}{\bi{2n}n}$$
and hence
\begin{equation}\label{k-1'}\sum_{k=1}^\infty\f{(-1)^k}{k\bi{2k}k^5}(1025k^5-2556k^4+2566k^3-1276k^2+321k-32)=-1.
\end{equation}
Observe that \eqref{k-1}$+668\times$\eqref{k-1'} gives
\begin{align*}&\sum_{k=1}^\infty\f{(-1)^k}{\bi{2k}k^5}(684700k^4-1418358k^3+1100639k^2-365392k+47327)
\\&\qquad=-\f{5467+64\zeta(3)}5,\end{align*}
which is equivalent to the desired (T9) since
$$47327-\f{5467}5=\f{231168}5.$$
This concludes our proof of Theorem \ref{Th1.3}. \qed

 \section{Other variants of Ramanujan-type series for $\f1{\pi}$}
 \setcounter{theorem}{0}
\setcounter{equation}{0}
\setcounter{conjecture}{0}

\begin{lemma} \label{223} Let $m\in\Z\sm\{0\}$ and $n\in\N$.

{\rm (i)} We have
\begin{align*}&\sum_{k=0}^n\f{\bi{2k}k^2\bi{3k}k}{m^k}\l((108-m)k^3+162k^2+78k+12\r)
\\&\qquad=6(2n+1)(3n+1)(3n+2)\f{\bi{2n}n^2\bi{3n}n}{m^n},
\end{align*}
\begin{align*}
&\sum_{k=0}^n\f{k\bi{2k}k^2\bi{3k}k}{m^k}\l((108-m)k^3+(162+m)k^2+78k+12\r)
\\&\qquad=6n(2n+1)(3n+1)(3n+2)\f{\bi{2n}n^2\bi{3n}n}{m^n},
\end{align*}
and
\begin{equation*}\begin{aligned}
&\sum_{k=0}^n\f{k^2\bi{2k}k^2\bi{3k}k}{m^k}\l((108-m)k^3+(162+2m)k^2+(78-m)k+12\r)
\\&\qquad\qquad=6n^2(2n+1)(3n+1)(3n+2)\f{\bi{2n}n^2\bi{3n}n}{m^n}.
\end{aligned}\end{equation*}

{\rm (ii)} We have
\begin{align*}&\sum_{k=0}^n\f{\bi{2k}k^2\bi{3k}k}{(k+1)^2m^k}
\l((108-m)k^3+(162-2m)k^2+(78-m)k+12\r)
\\&\quad=6(2n+1)(3n+1)(3n+2)\f{\bi{2n}n^2\bi{3n}n}{(n+1)^2m^n}
\end{align*}
and
\begin{align*}&\sum_{k=0}^n\f{k\bi{2k}k^2\bi{3k}k}{(k+1)^2m^k}
\l((108-m)k^3+(162-m)k^2+(78+m)k+m+12\r)
\\&\quad=6n(2n+1)(3n+1)(3n+2)\f{\bi{2n}n^2\bi{3n}n}{(n+1)^2m^n}.
\end{align*}
\end{lemma}
\Proof. It is easy to prove the five identities by induction on $n$.
We find them by the Gosper algorithm. \qed

\begin{theorem}\label{A-23} We have the following identities:
\begin{align*}
\sum_{k=0}^\infty\f{(150k^2-51k+89)k^3\bi{2k}k^2\bi{3k}k}{(-192)^k}&=\f{32\sqrt3}{15\pi},\tag{S4}
\\\sum_{k=0}^\infty\f{(1550k^2+3029k+ 1481)\bi{2k}{k+1}^2\bi{3k}k}{(-192)^k}&=-\f{80\sqrt3}{\pi},\tag{S4$'$}
\\\sum_{k=0}^\infty\f{(738k^2-4023k+745)k^3\bi{2k}k^2\bi{3k}k}{216^k}&=-\f{24\sqrt3}{\pi},\tag{S5}
\\\sum_{k=0}^\infty\f{(702k^2+1491k+787)\bi{2k}{k+1}^2\bi{3k}k}{216^k}&=\f{228\sqrt3}{\pi},\tag{S5$'$}
\\
\sum_{k=0}^\infty\f{(159426k^2-292761k+153995)k^3\bi{2k}k^2\bi{3k}k}{(-12)^{3k}}&=\f{96\sqrt3}{\pi},\tag{S6}
\\\sum_{k=0}^\infty\f{(604962k^2+1206195k+601183)\bi{2k}{k+1}^2\bi{3k}k}{(-12)^{3k}}&=-\f{6864\sqrt3}{\pi},\tag{S6$'$}
\end{align*}
\begin{align*}
\sum_{k=0}^\infty\f{(900k^2-2097k+929)k^3\bi{2k}k^2\bi{3k}k}{1458^k}&=-\f{27}{20\pi},\tag{S7}
\\\sum_{k=0}^\infty\f{(56025k^2+112584k+56551)\bi{2k}{k+1}^2\bi{3k}k}{1458^k}
&=\f{6615}{4\pi},\tag{S7$'$}
\end{align*}
\[
\sum_{k=0}^\infty\f{(173502k^2-354087k+181205)k^3\bi{2k}k^2\bi{3k}k}{(-8640)^k}=\f{32\sqrt{15}}{5\pi},
\tag{S8}\]
\[\begin{aligned}&\sum_{k=0}^\infty\f{(2086398k^2+4169997k+2083585)\bi{2k}{k+1}^2\bi{3k}k}{(-8640)^k}
=-\f{11504\sqrt{15}}{5\pi},\end{aligned}\tag{S8$'$}\]
\[\sum_{k=0}^\infty\f{(216711k^2-473742k+226495)k^3\bi{2k}k^2\bi{3k}k}{15^{3k}}=-\f{60\sqrt3}{\pi},
\tag{S9}\]
\[\begin{aligned}
&\sum_{k=0}^\infty\f{(620730k^2+1243839k+623087)\bi{2k}{k+1}^2\bi{3k}k}{15^{3k}}
=\f{16935\sqrt3}{4\pi},\end{aligned}\tag{S9$'$}\]
\[\sum_{k=0}^\infty\f{(36512550k^2-76893003k+39703217)k^3\bi{2k}k^2\bi{3k}k}{(-48)^{3k}}
=\f{768\sqrt3}{5\pi},\tag{S10}\]
\[\sum_{k=0}^\infty\f{Q_{-48}(k)\bi{2k}{k+1}^2\bi{3k}k}{(-48)^{3k}}
=-3538560\f{\sqrt3}{\pi},\tag{S10$'$}\]
where
$$Q_{-48}(k):=18009838350k^2+36017730237k+18007890593.$$
Also,
\[\begin{aligned}
&\sum_{k=0}^\infty\f{(30437550k^2-64346463k+33372157)k^3\bi{2k}k^2\bi{3k}k}{(-326592)^k}
=\f{864\sqrt7}{35\pi},\end{aligned}\tag{S11}\]
and
\[\begin{aligned}&\sum_{k=0}^\infty\f{Q(k)\bi{2k}{k+1}^2\bi{3k}k}{(-326592)^k}
=-11756880\f{\sqrt7}{7\pi},\end{aligned}\tag{S11$'$}\]
where
$$Q(k):=38518093350k^2+77034773577k+38516679853.$$
Moreover,
\[\begin{aligned}&\sum_{k=0}^\infty\f{(1227699410778k^-2613341265669k+1370770039375)k^3\bi{2k}k^2\bi{3k}k}{(-300)^{3k}}
\\&\qquad\qquad=12000\f{\sqrt3}{\pi},\end{aligned}\tag{S12}\]
and
\[\begin{aligned}&
\sum_{k=0}^\infty\f{Q_{-300}(k)\bi{2k}{k+1}^2\bi{3k}k}{(-300)^{3k}}
=-13499994000\f{\sqrt3}{\pi},\end{aligned}\tag{S12$'$}\]
where
\begin{align*}Q_{-300}(k)=&16746723121124538k^2+33493438799463639k
\\&+16746715678300691.
\end{align*}
\end{theorem}
\Proof. By Lemma \ref{223}(i), for any $m\in\Z$ with $|m|>108$, we have
\begin{align}\sum_{k=0}^\infty\f{\bi{2k}k^2\bi{3k}k}{m^k}\l((108-m)k^3+162k^2+78k+12\r)
&=0,
\\\sum_{k=0}^\infty\f{k\bi{2k}k^2\bi{3k}k}{m^k}\l((108-m)k^3+(162+m)k^2+78k+12\r)
&=0,
\end{align}
and
\begin{equation}
\sum_{k=0}^n\f{k^2\bi{2k}k^2\bi{3k}k}{m^k}\l((108-m)k^3+(162+2m)k^2+(78-m)k+12\r)
=0.\end{equation}

 We now prove (S4) in details. Taking $m=-192$ in the last three identities, we get
\begin{align}\label{192-0}\sum_{k=0}^\infty\f{\bi{2k}k^2\bi{3k}k}{(-192)^k}(50k^3+27k^2+13k+2)
=&0,
\\\label{192-1}\sum_{k=0}^\infty\f{k\bi{2k}k^2\bi{3k}k}{(-192)^k}(50k^3-5k^2+13k+2)
=&0,
\\\label{192-2}\sum_{k=0}^\infty\f{k^2\bi{2k}k^2\bi{3k}k}{(-192)^k}(50k^3-37k^2+45k+2)
=&0.
\end{align}
Combining the known identity
\begin{equation}\label{192o}\sum_{k=0}^\infty\f{\bi{2k}k^2\bi{3k}k}{(-192)^k}(5k+1)=\f{4}{\sqrt3\pi}
\end{equation}
with \eqref{192-0}, we get
$$\sum_{k=0}^\infty\f{k\bi{2k}k^2\bi{3k}k}{(-192)^k}(50k^2+27k+3)=-\f{8}{\sqrt3\pi}.$$
Combining this with \eqref{192-1}, we see that
$$\sum_{k=0}^\infty\f{k\bi{2k}k^2\bi{3k}k}{(-192)^k}(150k^3-(3\times5+2\times50)k^2+(3\times13-2\times27)k)
=\f{16\sqrt3}{\pi},$$
i.e.,
\begin{equation}\label{kk}\sum_{k=0}^\infty\f{k^2\bi{2k}k^2\bi{3k}k}{(-192)^k}
(30k^2-23k-3)=\f{16\sqrt3}{5\pi}.\end{equation}
\eqref{192-2}$\times 3$ plus \eqref{kk}$\times2$ yields
$$\sum_{k=0}^\infty\f{k^2\bi{2k}k^2\bi{3k}k}{(-192)^k}
(150k^3+(60-111)k^2+(135-46)k)=\f{32}{5\sqrt3\pi},$$
which is equivalent to (S4).

Next we prove (S4$'$). Applying \ref{223}(ii) with $m=-192$, we have
\begin{align}\label{0-192}\sum_{k=0}^\infty\f{\bi{2k}k^2\bi{3k}k}{(k+1)^2(-192)^k}
(50k^3+91k^2+45k+2)=&0,
\\\label{1-192}\sum_{k=0}^\infty\f{k\bi{2k}k^2\bi{3k}k}{(k+1)^2(-192)^k}
(50k^3+59k^2-19k-30)=&0.
\end{align}
As
$$50k^3+91k^2=45k+2-10(5k+1)(k+1)^2=-(19k^2+25k+8),$$
from \eqref{0-192} and \eqref{192o} we obtain
\begin{equation}\label{ab}\sum_{k=0}^\infty\f{\bi{2k}k^2\bi{3k}k}{(k+1)^2(-192)^k}(19k^2+25k+8)=\f{40\sqrt3}{3\pi}.
\end{equation}
Note that \eqref{0-192}$\times4-$\eqref{ab} yields
\begin{equation}\label{cc}\sum_{k=0}^\infty\f{k\bi{2k}k^2\bi{3k}k}{(k+1)^2(-192)^k}(40k^2+69k+31)=-\f{8\sqrt3}{3\pi}.
\end{equation}
Via \eqref{1-192}$\times31+30\times$\eqref{cc} we obtain
$$\sum_{k=0}^\infty\f{k^2\bi{2k}k^2\bi{3k}k}{(k+1)^2(-192)^k}(1550k^2+3029k+1481)=-\f{80\sqrt3}{\pi},$$
which is equivalent to the desired (S4$'$) since $\bi{2k}{k+1}=\f k{k+1}\bi{2k}k$.

All the remaining formulas in Theorem \ref{A-23} can be proved similarly. We omit the details. \qed
\qed

\begin{lemma} \label{Lem-24} Let $m\in\Z\sm\{0\}$ and $n\in\N$.

{\rm (i)} We have
\begin{align*}&\sum_{k=0}^n\f{\bi{2k}k^2\bi{4k}{2k}}{m^k}\l((256-m)k^3+384k^2+176k+24\r)
\\&\qquad=8(2n+1)(4n+1)(4n+3)\f{\bi{2n}n^2\bi{4n}{2n}}{m^n},
\end{align*}
\begin{align*}
&\sum_{k=0}^n\f{k\bi{2k}k^2\bi{4k}{2k}}{m^k}\l((256-m)k^3+(384+m)k^2+176k+24\r)
\\&\qquad=8n(2n+1)(4n+1)(4n+3)\f{\bi{2n}n^2\bi{4n}{2n}}{m^n},
\end{align*}
and
\begin{equation*}\begin{aligned}
&\sum_{k=0}^n\f{k^2\bi{2k}k^2\bi{4k}{2k}}{m^k}\l((256-m)k^3+(384+2m)k^2+(176-m)k+24\r)
\\&\qquad\qquad=8n^2(2n+1)(4n+1)(4n+3)\f{\bi{2n}n^2\bi{4n}{2n}}{m^n}.
\end{aligned}\end{equation*}

{\rm (ii)} We have
\begin{align*}&\sum_{k=0}^n\f{\bi{2k}k^2\bi{4k}{2k}}{(k+1)^2m^k}
\l((256-m)k^3+(384-2m)k^2+(176-m)k+24\r)
\\&\quad=8(2n+1)(4n+1)(4n+3)\f{\bi{2n}n^2\bi{4n}{2n}}{(n+1)^2m^n}
\end{align*}
and
\begin{align*}&\sum_{k=0}^n\f{k\bi{2k}k^2\bi{4k}{2k}}{(k+1)^2m^k}
\l((256-m)k^3+(384-m)k^2+(176+m)k+m+24\r)
\\&\quad=8n(2n+1)(4n+1)(4n+3)\f{\bi{2n}n^2\bi{4n}{2n}}{(n+1)^2m^n}.
\end{align*}
\end{lemma}
\Proof. It is easy to prove the five identities by induction on $n$.
We find them by the Gosper algorithm. \qed

Using Lemma \ref{Lem-24} and the method we prove Theorem \ref{A-23}, we obtain the following theorem.

\begin{theorem}\label{A24} We have the following identities:
\begin{align*}\sum_{k=0}^\infty\f{(854k^2-3639k+910)k^3\bi{2k}k^2\bi{4k}{2k}}{648^k}&=-\f{243}{14\pi},\tag{S13}
\\\sum_{k=0}^\infty\f{(2107k^2+4359k+2249)\bi{2k}{k+1}^2\bi{4k}{2k}}{648^k}&=\f{567}{\pi},\tag{S13$'$}
\\\sum_{k=0}^\infty\f{(400k^2-496k+327)k^3\bi{2k}k^2\bi{4k}{2k}}{(-1024)^k}&=\f9{5\pi},\tag{S14}
\\\sum_{k=0}^\infty\f{(19840k^2+39324k+19481)\bi{2k}{k+1}^2\bi{4k}{2k}}{(-1024)^k}&=-\f{1000}{\pi},\tag{S14$'$}
\end{align*}
\begin{align*}
\sum_{k=0}^\infty\f{(6016k^2-14856k+6215)k^3\bi{2k}k^2\bi{4k}{2k}}{48^{2k}}&=-\f{27\sqrt3}{4\pi},\tag{S15}
\\\sum_{k=0}^\infty\f{(75776k^2+152520k+76729)\bi{2k}{k+1}^2\bi{4k}{2k}}{48^{2k}}&
=\f{1746\sqrt3}{\pi},\tag{S15$'$}
\\
\sum_{k=0}^\infty\f{(1326650k^2-2420121k+1281559)k^3\bi{2k}k^2\bi{4k}{2k}}{(-63^2)^k}&=1944\f{\sqrt7}{5\pi},
\tag{S16}
\end{align*}
and
\[\begin{aligned}&\sum_{k=0}^\infty\f{(8242975k^2+16441878k+8198387)\bi{2k}{k+1}^2\bi{4k}{2k}}{(-63^2)^k}&
=-\f{106515\sqrt7}{2\pi}.\end{aligned}\tag{S16$'$}\]
We also have
\[\sum_{k=0}^\infty\f{(320k^2-678k+337)k^3\bi{2k}k^2\bi{4k}{2k}}{12^{4k}}=-\f{243\sqrt2}{10240\pi},\tag{S17}\]
\[\begin{aligned}&\sum_{k=0}^\infty\f{(6671360k^2+13350714k+6679321)\bi{2k}{k+1}^2\bi{4k}{2k}}{12^{4k}}
=\f{70065\sqrt2}{4\pi},\end{aligned}\tag{S17$'$}\]
\[\sum_{k=0}^\infty\f{(2576k^2-5136k+2635)k^3\bi{2k}k^2\bi{4k}{2k}}{(-3\times2^{12})^k}=\f{2\sqrt3}{7\pi},
\tag{S18}\]
\[\begin{aligned}
&\sum_{k=0}^\infty\f{(796544k^2+1591612k+795059)\bi{2k}{k+1}^2\bi{4k}{2k}}{(-3\times2^{12})^k}
=-\f{8176\sqrt3}{3\pi},\end{aligned}
\tag{S18$'$}\]
\[
\sum_{k=0}^\infty\f{(2428400k^2-5044368k+2584321)k^3\bi{2k}k^2\bi{4k}{2k}}{(-2^{10}3^4)^k}=\f{243}{5\pi},
\tag{S19}\]
\[\begin{aligned}&\sum_{k=0}^\infty\f{(2475740800k^2+4950772932k+ 2475031103)\bi{2k}{k+1}^2\bi{4k}{2k}}{(-2^{10}3^4)^k}
\\&\qquad\qquad=-\f{2238840}{\pi},\end{aligned}\tag{S19$'$}\]
\[\sum_{k=0}^\infty\f{(38400k^2-80696k+41609)k^3\bi{2k}k^2\bi{4k}{2k}}{28^{4k}}=-\f{49\sqrt3}{1080\pi},
\tag{S20}\]
\[\begin{aligned}&\sum_{k=0}^\infty\f{(1967513600k^2+3935104168k+ 1967590547)\bi{2k}{k+1}^2\bi{4k}{2k}}{28^{4k}}
\\&\qquad\qquad
=\f{3764915\sqrt3}{27\pi},\end{aligned}\tag{S20$'$}\]
\[\sum_{k=0}^\infty\f{(423089968k^2-891891888k +463383905)k^3\bi{2k}k^2\bi{4k}{2k}}{(-2^{14}3^45)^k}=\f{972\sqrt5}{35\pi},\tag{S21}\]
\[\sum_{k=0}^\infty\f{q(k)\bi{2k}{k+1}^2\bi{4k}{2k}}{(-2^{14}3^45)^k}
=-716633568\f{\sqrt5}{5\pi},\tag{S21$'$}\]
where
\begin{align*}q(k):=&28206829076608k^2+56413556154372k+28206727074305.
\end{align*}
We also have
\[\begin{aligned}&\sum_{k=0}^\infty\f{(643835623600k^2-1361740501968k+ 711617288021)k^3\bi{2k}k^2\bi{4k}{2k}}{(-2^{10}21^4)^k}
\\&\qquad\qquad=\f{11907}{5\pi},\end{aligned}\tag{S22}\]
and
\[\sum_{k=0}^\infty\f{q_k\bi{2k}{k+1}^2\bi{4k}{2k}}{(-2^{10}21^4)^k}=
-\f{263473491960}{\pi},\tag{S22$'$}\]
where $q_k$ denotes
\begin{align*}695911303499907200k^2+1391822523134211732k+695911219634175403
\end{align*}
Moreover,
\[\sum_{k=0}^\infty\f{(13462400k^2-28347528k+14689591)k^3\bi{2k}k^2\bi{4k}{2k}}{1584^{2k}}
=-\f{243\sqrt{11}}{140\pi},
\tag{S23}\]
and
\[\sum_{k=0}^\infty\f{u_k\bi{2k}{k+1}^2\bi{4k}{2k}}{1584^{2k}}
=16936290\f{\sqrt{11}}{\pi},\tag{S23$'$}\]
where
\begin{align*}u_k:=1868912998400k^2+3737843883096k+1868930883259.
\end{align*}
Also,
\[\begin{aligned}
\\&\sum_{k=0}^\infty\f{(234440315200k^2-497134511862k+260991361673)k^3\bi{2k}k^2\bi{4k}{2k}}{396^{4k}}
\\&\qquad\qquad=-\f{264627\sqrt2}{71680\pi},\end{aligned}\tag{S24}\]
and
\[\begin{aligned}&\sum_{k=0}^\infty\f{v_k\bi{2k}{k+1}^2\bi{4k}{2k}}{396^{4k}}
=90382094430705\f{\sqrt2}{4\pi},
\end{aligned}\tag{S24$'$}\]
where
\begin{align*}v_k:=&10422143035665206809600k^2+20844286081501973765862k
\\&+10422143045836766781623
\end{align*}
\end{theorem}

\begin{lemma} \label{Lem-36} Let $m\in\Z\sm\{0\}$ and $n\in\N$.

{\rm (i)} We have
\begin{align*}&\sum_{k=0}^n\f{\bi{2k}k\bi{3k}k\bi{6k}{3k}}{m^k}\l((1728-m)k^3+2592k^2+1104k+120\r)
\\&\qquad=24(2n+1)(6n+1)(6n+5)\f{\bi{2n}n\bi{3n}n\bi{6n}{3n}}{m^n},
\end{align*}
\begin{align*}
&\sum_{k=0}^n\f{k\bi{2k}k\bi{3k}k\bi{6k}{3k}}{m^k}\l((1728-m)k^3+(2592+m)k^2+1104k+120\r)
\\&\qquad=24n(2n+1)(6n+1)(6n+5)\f{\bi{2n}n\bi{3n}n\bi{6n}{3n}}{m^n},
\end{align*}
and
\begin{equation*}\begin{aligned}
&\sum_{k=0}^n\f{k^2\bi{2k}k\bi{3k}k\bi{6k}{3k}}{m^k}\l((1728-m)k^3+(2592+2m)k^2+(1104-m)k+120\r)
\\&\qquad\qquad=24n^2(2n+1)(6n+1)(6n+5)\f{\bi{2n}n\bi{3n}n\bi{6n}{3n}}{m^n}.
\end{aligned}\end{equation*}

{\rm (ii)} We have
\begin{align*}&\sum_{k=0}^n\f{\bi{2k}k\bi{3k}k\bi{6k}{3k}}{(k+1)^2m^k}
\l((1728-m)k^3+(2592-2m)k^2+(1104-m)k+120\r)
\\&\quad=24(2n+1)(6n+1)(6n+5)\f{\bi{2n}n\bi{3n}n\bi{6n}{3n}}{(n+1)^2m^n}
\end{align*}
and
\begin{align*}&\sum_{k=0}^n\f{k\bi{2k}k\bi{3k}k\bi{6k}{3k}}{(k+1)^2m^k}
\l((1728-m)k^3+(2592-m)k^2+(1104+m)k+m+120\r)
\\&\quad=24n(2n+1)(6n+1)(6n+5)\f{\bi{2n}n^2\bi{4n}{2n}}{(n+1)^2m^n}.
\end{align*}
\end{lemma}
\Proof. It is easy to prove the five identities by induction on $n$.
We find them by the Gosper algorithm. \qed

For any $k\in\N$, we clearly have
$$\bi{3k}{k+1}=\f{2k}{k+1}\bi{3k}k$$
which is similar to $\bi{2k}{k+1}=\f k{k+1}\bi{2k}k$.
In view of this and Lemma \ref{Lem-36}, by our method to prove Theorem \ref{A-23} we can establish the following theorem.

\begin{theorem}\label{A236} {\rm (i)} We have the following identities:
\begin{align*}
\sum_{k=0}^\infty\f{(23968k^2-71188k+25539)k^3\bi{2k}k\bi{3k}k\bi{6k}{3k}}{20^{3k}}&=-\f{1875\sqrt5}{56\pi},
\tag{S25}
\\\sum_{k=0}^\infty\f{(269696k^2+545164k+275383)\bi{2k}{k+1}\bi{3k}{k+1}\bi{6k}{3k}}{20^{3k}}
&=\f{50750\sqrt5}{3\pi},\tag{S25$'$}
\\\sum_{k=0}^\infty\f{(90706k^2-168589k+88872)k^3\bi{2k}k\bi{3k}k\bi{6k}{3k}}{(-2^{15})^k}
&=\f{150\sqrt2}{7\pi},\tag{S26}
\end{align*}
and
\[\sum_{k=0}^\infty\f{(5821200k^2+11623266k+5801911)\bi{2k}{k+1}\bi{3k}{k+1}\bi{6k}{3k}}{(-2^{15})^k}
=-\f{261184\sqrt2}{3\pi}.\tag{S26$'$}\]

{\rm (ii)} We have
\[\sum_{k=0}^\infty\f{(457326k^2-308241k+342923)k^3\bi{2k}k\bi{3k}k\bi{6k}{3k}}{(-15)^{3k}}
=\f{5000\sqrt{15}}{7\pi},\tag{S27}\]
and
\[\sum_{k=0}^\infty\f{(933849k^2+1843866k+910277)\bi{2k}{k+1}\bi{3k}{k+1}\bi{6k}{3k}}{(-15)^{3k}}=-\f{27125\sqrt{15}}{\pi},
\tag{S27$'$}\]
Moreover,
\[\sum_{k=0}^\infty\f{(12826k^2-27741k+13298)k^3\bi{2k}k\bi{3k}k\bi{6k}{3k}}{(2\times30^3)^k}
=-\f{625\sqrt{15}}{972\pi},\tag{S28}\]
\[\begin{aligned}&\sum_{k=0}^\infty\f{(4941882k^2+9895613k+4953661)\bi{2k}{k+1}\bi{3k}{k+1}\bi{6k}{3k}}{(2\times30^3)^k}
\\&\qquad\qquad=\f{56375\sqrt{15}}{3\pi}.\end{aligned}\tag{S28$'$}\]

{\rm (iii)} We have
\[\sum_{k=0}^\infty\f{(18126342k^2-37421775k+19111480)k^3\bi{2k}k\bi{3k}k\bi{6k}{3k}}{(-96)^{3k}}
=\f{50\sqrt6}{\pi},\tag{S29}\]
and
\[\begin{aligned}&\sum_{k=0}^\infty\f{(6802059888k^2+13603203918k+6801143345)\bi{2k}{k+1}\bi{3k}{k+1}\bi{6k}{3k}}{(-96)^{3k}}
\\&\qquad\quad=-2358976\f{\sqrt6}{\pi}.\end{aligned}
\tag{S29$'$}\]

 {\rm (iv)} We have
\[\sum_{k=0}^\infty\f{(2248722k^2-4689621k+2357878)k^3\bi{2k}k\bi{3k}k\bi{6k}{3k}}{66^{3k}}
=-\f{275\sqrt{33}}{28\pi},\tag{S30}\]
 \[\begin{aligned}
&\sum_{k=0}^\infty\f{(71477721k^2+142985853k+71508077)\bi{2k}{k+1}\bi{3k}{k+1}\bi{6k}{3k}}{66^{3k}}
=32956\f{\sqrt{33}}{\pi},\end{aligned}\tag{S30$'$}
\]
\[\sum_{k=0}^\infty\f{(2161071858k^2-4497745053k+2312761384)k^3\bi{2k}k\bi{3k}k\bi{6k}{3k}}
{(-3\times160^3)^k}=\f{1250\sqrt{30}}{9\pi},\tag{S31}\]
and
\[\sum_{k=0}^\infty\f{w_k\bi{2k}{k+1}\bi{3k}{k+1}\bi{6k}{3k}}{(-3\times160^3)^k}
=-273064000\f{\sqrt{30}}{3\pi},\tag{S31$'$}\]
where
\begin{align*}w_k:=&8126882714192k^2+16253686120778k+8126803400291
\end{align*}

{\rm (v)} We have
\[\begin{aligned}
\\&\sum_{k=0}^\infty\f{(273732850062k^2-572136425667k +296241014776)k^3\bi{2k}k\bi{3k}k\bi{6k}{3k}}{(-960)^{3k}}
\\&\qquad\qquad=\f{5000\sqrt{15}}{21\pi},\end{aligned}\tag{S32}\]
and
\[\sum_{k=0}^\infty\f{x_k\bi{2k}{k+1}\bi{3k}{k+1}\bi{6k}{3k}}{(-960)^{3k}}
=-235929568000\f{\sqrt{15}}{3\pi},\tag{S32$'$}\]
where
\begin{align*}x_k:=&357405080886027216k^2+ 714810113296569594k+357405032410460843.
\end{align*}
Also,
\[\sum_{k=0}^\infty\f{(6491502k^2-13521457k+6955771)k^3\bi{2k}k\bi{3k}k\bi{6k}{3k}}{255^{3k}}
=-\f{42500\sqrt{255}}{413343\pi},\tag{S33}\]
and
\[
\sum_{k=0}^\infty\f{y_k
\bi{2k}{k+1}\bi{3k}{k+1}\bi{6k}{3k}}{255^{3k}}
=2349045125\f{\sqrt{255}}{54\pi}.\tag{S33$'$}\]
where
\begin{align*}y_k:=&15274005325299k^2+30548121249166k+15274115917127.
\end{align*}

{\rm (vi)}
We have
\[\sum_{k=0}^\infty\f{p(k)k^3\bi{2k}k\bi{3k}k\bi{6k}{3k}}{(-5280)^{3k}}
=\f{13750\sqrt{330}}{7\pi}\tag{S34}\]
and
\[
\begin{aligned}
&\sum_{k=0}^\infty\f{q(k)\bi{2k}{k+1}\bi{3k}{k+1}\bi{6k}{3k}}{(-5280)^{3k}}
=-107945164712000\f{\sqrt{330}}{\pi},
\end{aligned}\tag{S34$'$}\]
where
\begin{align*}
p(k)=&2417841826472898k^2-5065781116806693k +2634006768739304
\end{align*}
and
\begin{align*}
q(k)=&382825470402996808454064k^2+765650940493903329223926k
\\&+382825470090906516313597.
\end{align*}
Also,
\[
\begin{aligned}
&\sum_{k=0}^\infty\f{P(k)k^3\bi{2k}k\bi{3k}k\bi{6k}{3k}}{(-640320)^{3k}}
=-\f{83\sqrt{10005}}{\pi}
\end{aligned}\tag{S35}\]
and
\[\begin{aligned}
\sum_{k=0}^\infty\f{Q(k)\bi{2k}{k+1}\bi{3k}{k+1}\bi{6k}{3k}}{(-640320)^{3k}}
=-46696654461704580256000\f{\sqrt{10005}}{\pi},
\end{aligned}\tag{S35$'$}\]
where
\begin{align*}P(k):=&3726784819871553194540063287782k^2
\\&-7783860761103294083667021327391k
\\&+4057075941237594195269253626425
\end{align*}
and
\begin{align*}Q(k):=&1626388893999999577578620229159002547888k^2 \\&+3252777787999998411771125249563926009942k
\\&+1626388893999998834192505020394081197549.
\end{align*}
\end{theorem}

\section{Proof of Theorem \ref{Th1.4}}
 \setcounter{theorem}{0}
\setcounter{equation}{0}
\setcounter{conjecture}{0}

\medskip
\noindent{\it Proof of Theorem} 1.4(i). Let $u_n=\Domb(n)/64^n$
and $v_n=\Domb(n)/(-32)^n$ for $n\in\N$. By the Zeilberger algorithm,
we find the recurrence
$$(n+1)^3u_n-2(2n+3)(5n^2+15n+12)u_{n+1}+64(n+2)^3u_{n+2}=0$$
and
$$(n+1)^3v_n+(2n+3)(5n^2+15n+12)v_{n+1}+16(n+2)^3v_{n+2}=0.$$
Thus,
\begin{align*}0=&\sum_{n=0}^\infty(n+1)^3u_n-2\sum_{n=0}^\infty(2(n+1)+1)(5(n+1)^2+5(n+1)+2)u_{n+1}
\\&+64\sum_{n=0}^\infty(n+2)^3u_{n+2}
\\=&\sum_{k=0}^\infty(k+1)^3u_k-2\sum_{k=1}^\infty(2k+1)(5k^2+5k+2)u_k+64\sum_{k=2}^\infty k^3u_k
\\=&\sum_{k=0}^\infty\l((k+1)^3-2(2k+1)(5k^2+5k+2)+64k^3\r)u_k
\\&+2(2\times0+1)(5\times0^2+5\times0+2)u_0-64(0^3u_0+1^3u_1)
\\=&3\sum_{k=0}^\infty(15k^3-9k^2-5k-1)u_k
\end{align*}
and hence
\begin{equation}\label{u-k}\sum_{k=0}^\infty(15k^3-9k^2-5k-1)u_k=0.
\end{equation}
Similarly,
\begin{align*}0=&\sum_{k=0}^\infty\l((k+1)^3+(2k+1)(5k^2+5k+2)+16k^3\r)v_k
\\=&3\sum_{k=0}^\infty(9k^3+6k^2+4k+1)v_k
\end{align*}
and hence
\begin{equation}\label{v-k}\sum_{k=0}^\infty(9k^3+6k^2+4k+1)v_k=0.
\end{equation}
Combining \eqref{u-k} and \eqref{v-k} with the two known identities in \eqref{64-32},
we immediately get the desired \eqref{D64} and \eqref{D-32}. \qed

\medskip
\noindent{\it Proof of Theorem} 1.4(ii).
 The second part of Theorem \eqref{Th1.4} can be proved by the same method we prove
Theorem \eqref{Th1.4}(i). Here we just provide a proof of \eqref{Y36} in details.
Let $w_n=f_n^{(4)}/36^n$ for $n\in\N$. By the Zeilberger algorithm, we find the recurrence
$$324(n+2)^3w_{n+2}-18(2n+3)(3n^2+9n+7)w_{n+1}-(n+1)(4n+3)(4n+5)w_n=0.$$
Thus,
\begin{align*}0=&324\sum_{n=0}^\infty(n+2)^3w_{n+2}-18\sum_{n=0}^\infty(2(n+1)+1)(3(n+1)^2+3(n+1)+1)w_{n+1}
\\&-\sum_{n=0}^\infty(n+1)(4n+3)(4n+5)w_n
\\=&324\sum_{k=2}^\infty k^3w_k-18\sum_{k=1}^\infty(2k+1)(3k^2+3k+1)w_k
\\&-\sum_{k=0}^\infty(k+1)(4k+3)(4k+5)w_k
\\=&\sum_{k=0}^\infty\l(324k^3-18(2k+1)(3k^2+3k+1)-(k+1)(4k+3)(4k+5)\r)w_k
\\&-324\times1^3w_1+18(2\times0+1)(3\times0^2+3\times0+1)
\end{align*}
and hence
\begin{equation}\label{200}\sum_{k=0}^\infty(200k^3-210k^2-137k-33)w_k=0.
\end{equation}
By \eqref{Yang},
\begin{equation}\label{Yang'}\sum_{k=0}^\infty 33(4k+1)w_k=\f{33\times18}{\sqrt{15}\,\pi}.
\end{equation}
Adding \eqref{200} and \eqref{Yang'}, we immediately obtain \eqref{Y36}. \qed
\medskip

Our method for proving parts (i)-(ii) of Theorem 1.4 also works for parts (iii)-(iv)
of Theorem 1.4. We omit the proof details.

\section{Conjectures for series only involving binomial coefficients}

Motivated by the author's conjectural identity
$$\sum_{k=1}^\infty\f{48^k}{k(2k-1)\bi{2k}k\bi{4k}{2k}}=\f {15}2K$$
(cf. \cite{S14d,S-11}), we make the following conjecture.

\begin{conjecture} We have
$$\sum_{k=0}^\infty\f{(4k^2-30k+17)48^k}{(2k-1)\bi{2k}k\bi{4k}{2k}}=\f{45K-44}2.$$
\end{conjecture}
\begin{remark} By induction, for each $n=1,2,3,\ldots$ we have
$$\sum_{k=1}^n\f{(4k^2-28k+3)48^k}{k\bi{2k}k\bi{4k}{2k}}=12-\f{12(n+1)48^n}{\bi{2n}n\bi{4n}{2n}}$$
and
$$\sum_{k=1}^n\f{(4k^2-40k-9)48^k}{\bi{2k}k\bi{4k}{2k}}=-12+\f{12(n+1)^248^n}{\bi{2n}n\bi{4n}{2n}}.
$$
\end{remark}

In view of Theorem \ref{Th1.3} and the author's conjectural identity
\eqref{Sun-zeta(3)}, we pose the following conjecture obtained via the PSLQ algorithm.

\begin{conjecture} \label{conj-1} {\rm (i)} We have
\begin{align*}\sum_{k=1}^\infty\f{(-1)^{k-1}}{\bi{2k}{k+1}^5}P(k)=\f{8103654862170335619+ 4368545100830839178\zeta(3)}5
\end{align*}
where
\begin{align*}P(k)=&54430524632163842275k^4-132483674356197881281k^3
\\&+121816306858962351437k^2-47590274284953796032k
\\&+6700636215039814272.
\end{align*}

{\rm (ii)} We have
\[\begin{aligned}&\sum_{k=0}^\infty\f{(-64)^k}{\bi{2k}k^4\bi{3k}k}
(676704k^4-1205388k^3+1140374k^2-152237k+78797)
\\&\qquad=35733-1344\zeta(3),\end{aligned}\]
$$\sum_{k=1}^\infty\f{(-64)^{k-1}}{\bi{2k}{k+1}^4\bi{3k}k}Q(k)=
7(14488697756+4718909979\zeta(3))$$
and
$$\sum_{k=1}^\infty\f{(-64)^kR(k)}{\bi{2k}{k+1}^4\bi{3k}{k+1}}
=-86166921288937568-74477398755902744\zeta(3),$$
where
\begin{align*}Q(k)=&9152858507744k^4 -18103487906940k^3+ 16104889340010k^2
\\&-5519172201903k+ 668801335410
\end{align*}
and
\begin{align*}R(k)=&152571345867547488k^4-325445013351260332k^3
\\&+295511129648313866k^2-106449469340961699k
\\&+13378286508841890
\end{align*}
\end{conjecture}
\begin{remark} In the spirit of our method to prove Theorems 1.1-1.3
as well as the general algorithm to deduce series of type $S$,  parts (i) and (ii) of Conjecture \ref{conj-1}
look equivalent to \eqref{AZ-zeta(3)} and \eqref{Sun-zeta(3)} respectively
 but we have not done this in details.
\end{remark}

Motivated by B. Gourevich's conjectural identity \eqref{Go} and Guillera's conjectural identity
$$\sum_{k=1}^\infty\f{256^k}{k^7\bi{2k}k^7}(21k^3-22k^2+8k-1)=\f{\pi^4}8$$
(cf. \cite[Section 4]{Gu-Exp}), we pose the following conjecture via the PSLQ algorithm.

\begin{conjecture} {\rm (i)} We have
\begin{align*}
\sum_{k=1}^\infty \f{S(k)}{2^{20k}}k^7\bi{2k}k^7=&-\f1{24\pi^3}
\end{align*}
and
\begin{align*}\f1{32}\sum_{k=1}^\infty\f{T(k)}{2^{20k}}\bi{2k}{k+1}^7
=&75570231394467396545747200
\\&-\f{2343145585133056805845704703}{\pi^3}
\end{align*}
where
\begin{align*}S(k):=&56448k^6-347200k^5+854280k^4
\\&-1145956k^3+851214k^2-339967k+56160
\end{align*}
and
\begin{align*}T(k)=&362334901725047543340617856k^6+1557477795272579082461315904k^5
        \\&+2387296377854823511932510816k^4+1269650426797215833274563064k^3
        \\&-372612359665735835469802516k^2-620395969622808879309367722k
        \\&-170581863683533821644571735
\end{align*}

{\rm (ii)} We have
\begin{align*}&\sum_{k=0}^\infty\f{U(k)256^k}{\bi{2k}{k}^7}=\f{9984-\pi^4}{48}
\end{align*}
and
\begin{align*}\sum_{k=1}^\infty\f{V(k)256^k}{\bi{2k}{k+1}^7}
=\f{10303458455020089696294336-47118227214655104697325\pi^4}{48},
\end{align*}
where
\begin{align*}U(k):=14112k^6-44464k^5+49490k^4-41069k^3+8155k^2-4749k-210.
\end{align*}
and
\begin{align*}V(k):=&105195631551721406964324k^6-375522390327670972174376k^5
\\&+508010030024769047270138k^4-376989573186149736346723k^3
\\&+156139327775481582503524k^2-33942569977747809706722k
 \\&+3027100288061502033993.
 \end{align*}
\end{conjecture}

\section{New type series involving generalized central trinomial coefficients}

For $b,c\in\Z$ and $n\in\N$, the generalized central trinomial coefficient
$T_n(b,c)$ denotes the coefficient of $x^n$ in the expansion of $(x^2+bx+c)^n$.
The author \cite{S14c,S-11,S20,S20p} posed totally 9 types of conjectural series for $1/\pi$ involving
generalized central trinomial coefficients. Here we consider their variants of type $S$.

\begin{conjecture} We have
\[\sum_{k=1}^\infty\f{3054600k^2-16826114k+11236485}{(-256)^k}
        k^3\bi{2k}k^2T_k(1,16)=-\f{1952307}{5\pi},\tag{I1S}\]
\[\sum_{k=1}^\infty\f{357600k^2-239434k+401075}{(-1024)^k}
        k^3\bi{2k}k^2T_k(34,1)=\f{25983}{10\pi},\tag{I2S}\]
\[\begin{aligned}&\sum_{k=1}^\infty\f{28823880k^2-740215234k+3516311133}{4096^k}k^3\bi{2k}k^2T_k(194,1)
        \\&\qquad\qquad=-\f{152854918}{3\pi},\end{aligned}\tag{I3S}\]
\[\begin{aligned}&\sum_{k=1}^\infty\f{1336776k^2-5896258k+4457117}{4096^k}
        k^3\bi{2k}k^2T_k(62,1)
       =-\f{188698\sqrt3}{7\pi}.\end{aligned}\tag{I4S}\]
Also,
\[\begin{aligned}&\sum_{k=1}^\infty\f{12k^2-28k-9}{256^k}
        k\bi{2k}k^2T_k(8,-2)=-\f{4(\sqrt{8+6\sqrt2}+3\root4\of{2})}{3\pi},
   \end{aligned}\]
\[\begin{aligned}&\sum_{k=1}^\infty\f{162k^2-969k+872}{256^k}
        k^2\bi{2k}k^2T_k(8,-2)=\f{708\root4\of{2}-1951\sqrt{8+6\sqrt2}}{12\pi},
   \end{aligned}\]
and
\[\begin{aligned}&\sum_{k=1}^\infty\f{31392k^2-277274k+594637}{256^k}
        k^3\bi{2k}k^2T_k(8,-2)
        \\&\quad=\f{90563\sqrt{8+6\sqrt2}-221844\,\root4\of{2}}{12\pi}.
   \end{aligned}\tag{I5S}\]
\end{conjecture}

\begin{conjecture} We have the following identities:

\[\begin{aligned}&\sum_{k=1}^\infty\f{k^3(3401775k^2-37884933k+65097406)}{972^k}
\bi{2k}k\bi{3k}kT_k(18,6)
\\&\qquad\quad=\f{158875\sqrt3}{\pi},\end{aligned}
\tag{II1S}\]
\[\begin{aligned}&\sum_{k=1}^\infty\f{k^3a_k}{1000^k}
\bi{2k}k\bi{3k}kT_k(10,1)=-\f{3441159450\sqrt3}{7\pi}\end{aligned}
\tag{II2S}\]
where
$$a_k=34073404820k^2-166944861551k+136909066683,$$
\[\begin{aligned}&\sum_{k=1}^\infty\f{k^3a_k'}{18^{3k}}
\bi{2k}k\bi{3k}kT_k(198,1)=-\f{1285805325750\sqrt3}{\pi}\end{aligned}
\tag{II3S}\]
where
$$a_k'=138317121900k^2-11341624063599k+169057291391203,$$
\[\begin{aligned}&\sum_{k=1}^\infty\f{515565k^2-1452888k+614707}{24^{3k}}
k^3\bi{2k}k\bi{3k}kT_k(26,729)\\&\qquad\quad=-\f{64}{98415\pi}(195423\sqrt3+532925\sqrt{15}),
\end{aligned}\tag{II10S}\]
\[\begin{aligned}&\sum_{k=1}^\infty\f{1326294k^2-4598217k+2285731}{(-5400)^k}
k^3\bi{2k}k\bi{3k}kT_k(70,3645)\\&\qquad\quad=-\f{2}{2187\pi}(78460655\sqrt3 +13402757\sqrt{15}),
\end{aligned}\tag{II11S}\]
\[\begin{aligned}&\sum_{k=1}^\infty\f{4076131815k^2-7828831071k+7124292568}{(-13500)^k}
k^3\bi{2k}k\bi{3k}kT_k(40,1548)\\&\qquad\quad=\f{25}{91854\pi}(32473732248\sqrt3-33657641611\sqrt6)
\end{aligned}\tag{II12S}\]

\[\begin{aligned}&\sum_{k=1}^\infty\f{17276571k^2-9174528k+19362029}{(-675)^k}
k^3\bi{2k}k\bi{3k}kT_k(15,-5)\\&\qquad\quad=\f{159680563\sqrt{15}}{3456\pi},
\end{aligned}\tag{II13S}\]
\[\begin{aligned}&\sum_{k=1}^\infty\f{5837598k^2-6981399k+6325061}{(-1944)^k}
k^3\bi{2k}k\bi{3k}kT_k(18,-3)\\&\qquad\quad=\f{11336\sqrt{3}}{\pi},
\end{aligned}\tag{II14S}\]
\end{conjecture}

\begin{conjecture} We have the following identities:
\begin{align*}&\sum_{k=1}^\infty\f{k(133120k^2-28704k+2669)}{(-168^2)^k}\bi{2k}k\bi{4k}{2k}T_k(7,4096)
\\&\qquad\quad=\f{21(14481\sqrt{42}+5650\sqrt{210})}{256\pi},
\end{align*}
\begin{align*}&\sum_{k=1}^\infty\f{k^2(2733056k^2-1413552k+1735253)}{(-168^2)^k}\bi{2k}k\bi{4k}{2k}T_k(7,4096)
\\&\qquad=\f{63(26395881\sqrt{42}+13904950\sqrt{210})}{40960\pi},
\end{align*}
and
\[\begin{aligned}&\sum_{k=1}^\infty\f{k^3b_k}{(-168^2)^k}\bi{2k}k\bi{4k}{2k}T_k(7,4096)
\\=&\f{189}{1310720\pi}(181371717913\sqrt{42}-200347079650\sqrt{210}),
\end{aligned}\tag{III5S}\]
where
$$b_k=17768990720k^2-81509977344k+73509891901.$$
\end{conjecture}
\begin{remark} We omit many other conjectural identities of type $S$ arising from the author's series for $1/\pi$ of type III (cf. \cite{S14c}).
\end{remark}

\begin{conjecture} We have the following identities:
\[\sum_{k=1}^\infty\f{k^3c_k}{(-48^2)^k}\bi{2k}k^2T_{2k}(7,1)=
-\f{7110645399}{5\pi}\tag{IV1S}
\]
and
\[\begin{aligned}&\sum_{k=1}^\infty\f{k^3d_k}{4608^k}\bi{2k}k^2T_{2k}(10,-2)
=-\f{20525967\sqrt{6}}{4\pi},
\end{aligned}\tag{IV19S}\]
where $$c_k=253607171350k^2-231835223289k+295912225514$$
and
$$d_k=220693694k^2-1203431385k+1673732470.$$
\end{conjecture}
\begin{remark} We omit many other conjectural identities of type $S$ arising from the author's series for $1/\pi$ of type IV (cf. \cite{S14c}).
\end{remark}

\begin{conjecture} We have
\[\sum_{k=1}^\infty\f{ke_k^{(1)}}{(-240)^{3k}}\bi{2k}k\bi{3k}kT_{3k}(62,1)
=-\f{318713456\sqrt{105}}{49\pi}\]
where
$$e_k^{(1)}:=24501204k^2+511623432k+17450401,$$
\[\sum_{k=1}^\infty\f{k^2e_k^{(2)}}{(-240)^{3k}}\bi{2k}k\bi{3k}kT_{3k}(62,1)
=-\f{11867637235544\sqrt{105}}{343\pi}\]
where
$$e_k^{(2)}:=1543522869252k^2+354135873240k+3269250841643,$$
and
\[\sum_{k=1}^\infty\f{k^3e_k^{(3)}}{(-240)^{3k}}\bi{2k}k\bi{3k}kT_{3k}(62,1)
=-\f{127447736489940892\sqrt{105}}{2401\pi}
\tag{V1S}\]
where
$$e_k^{(3)}:=41310253635000948k^2-37929668348178936k+56105265685716343.$$
Also,
\[\sum_{k=1}^\infty\f{(k^2+1)e_k}{(-240)^{3k}}\bi{2k}k\bi{3k}kT_{3k}(62,1)
=-\f{429939571410644\sqrt{105}}{343\pi},\]
where
$$e_k:=263893057256556k^2+46715038856064k+6234077522125.$$
\end{conjecture}
\begin{remark} This is motivated by the author's conjectural identity
\[\sum_{k=0}^\infty\f{1638k+277}{(-240)^{3k}}\bi{2k}k\bi{3k}kT_{3k}(62,1)=\f{44\sqrt{105}}{\pi}\tag{V1}\]
(cf. \cite{S14c,S-11}) confirmed in \cite{WZ}. We also conjecture that
$$\sum_{k=0}^\infty(47675628k^3+995541624k^2+437210809k+79678364)\f{\bi{2k}k\bi{3k}kT_{3k}(62,1)}{(-240)^{3k}}=0.$$
\end{remark}

\begin{conjecture} {\rm (i)} We have
\[\begin{aligned}&\sum_{k=1}^\infty\f{5939142726k^2-50217677843k+47679239989}{(2^{11}3^3)^k}k^3T_k(10,11^2)^3
\\&\qquad\qquad=\f{125366398162515\sqrt2}{644204\pi},\end{aligned}
\tag{VI1S}\]

{\rm (ii)} We have
\[\begin{aligned}&\sum_{k=1}^\infty\f{k^3f(k)}{(-80)^{3k}}T_k(22,21^2)^3
=-\f{317993969514116005\sqrt5}{64827\pi},\end{aligned}
\tag{VI2S}\]
where
\begin{align*}f(k):=&182075646906594k^2-658193121766971k+498776294291290.
\end{align*}

{\rm (iii)} We have
\[\begin{aligned}&\sum_{k=1}^\infty\f{k^3g(k)}{(-288)^{3k}}T_k(62,95^2)^3
=-\f{9}{30008125\pi}(a\sqrt2
        +b\sqrt{14}),\end{aligned}
\tag{VI3S}\]
where
$$a=2778571005723952224834213458\ ,\ b=2335639592156477133727790625,$$
and
\begin{align*}g(k):=&29814661986490996566930k^2-138270524135932079678425k
\\&+91756685418870140080439.
\end{align*}
\end{conjecture}

\begin{conjecture} {\rm (i)} We have
\[\begin{aligned}&\sum_{k=1}^\infty\f{k^3h_k}{450^{k}}\bi{2k}kT_k(6,2)^2
=-\f{1525516918799600750400}{7\pi},\end{aligned}
\tag{VII1S}\]
where
\begin{align*}h_k:=&1717562453635471595698k^2-20741022469043431508721k
\\&+29930651775020896516895.
\end{align*}
Also,
\[\begin{aligned}&\sum_{k=1}^\infty\f{761948702208k^2-4717770389584k+2473620838841}{28^{2k}}k^3\bi{2k}kT_k(4,9)^2
\\=&-\f{49}{46656\pi}(1744223791168\sqrt3+307854304805\sqrt6),
\end{aligned}\tag{VII2S}\]
and
\[\begin{aligned}&\sum_{k=1}^\infty\f{k^3\ell_k}{22^{2k}}\bi{2k}kT_k(5,1)^2
=-\f{861143519145937597955\sqrt7}{1536\pi},\end{aligned}
\tag{VII3S}\]
where
\begin{align*}\ell_k:=&25818737554793894400k^2-148156855009332208624k
\\&+136699546718553502681.
\end{align*}

{\rm (ii)} We have
\[\begin{aligned}&\sum_{k=1}^\infty\f{k^3m_k}{46^{2k}}\bi{2k}kT_k(7,1)^2
=-\f{675829671880947233751995193\sqrt7}{2458624\pi},\end{aligned}
\tag{VII4S}\]
where
\begin{align*}m_k:=&60405940502779625588736k^2-194226790034061510982928k
\\&+113445302068806184048447.
\end{align*}
Also,
\[\begin{aligned}&\sum_{k=1}^\infty\f{23575340288k^2-30724817208k+70092391411}{(-108)^k}k^3\bi{2k}kT_k(3,-3)^2
\\&\qquad\quad=\f{83673484461\sqrt7}{512\pi},\end{aligned}
\tag{VII5S}\]
and
\[\begin{aligned}&\sum_{k=1}^\infty\f{k^3n_k}{(-5177196)^k}\bi{2k}kT_k(171,-171)^2
\\=&\f{3506889487729535368666818898433157048375\sqrt7}{1024\pi},
\end{aligned}\tag{VII6S}\]
where
\begin{align*}n_k:=&10084570997509032944421348793965204251648k^2
\\&-19349061584246199246547660635419016030648k
\\&+10272620865354618637287798800250376143247.
\end{align*}
Moreover,
\[\begin{aligned}&\sum_{k=1}^\infty\f{k^3p(k)}{434^{2k}}\bi{2k}kT_k(73,576)^2
\\=&-\f{1941824970252133568742562638158781654482433600879227\sqrt6}{28311552\pi},
\end{aligned}\tag{VII7S}\]
where
\begin{align*}p(k):=&8486970607342072180410168982128577900986728448k^2
\\&-33680299486772296440807451280403962121626343552k
\\&+17329177465492357868335964242267514081171615873.
\end{align*}
\end{conjecture}

\begin{conjecture} We have
\[\begin{aligned}&\sum_{k=1}^\infty\f{4158282880k^2-22813966696k+40778481375}{(-50)^k}
        k^3T_k(4,1)T_k(1,-1)^2
        \\&\quad=-\f{29231055627965\sqrt{15}}{124416\pi},
\end{aligned}\tag{VIII1S}\]
\[\begin{aligned}&\sum_{k=1}^\infty\f{q(k)}{3240^k}
        k^3T_k(7,1)T_k(10,10)^2
        =-\f{2602848822349382206432377936\sqrt{5}}{7\pi},
\end{aligned}\tag{VIII2S}\]
where
\begin{align*}q(k):=&498872838526220266559278065k^2-7594176317285717732725105481k
        \\&+15117694789961557184766131894.
\end{align*}
Also,
\[\sum_{k=1}^\infty\f{r_k}{(-2430)^k}k^3T_k(8,1)T_k(5,-5)^2=\f{4114638030580119213489\sqrt{15}}{7168\pi},
\tag{VIII3S}\]
and
\[\begin{aligned}&\sum_{k=1}^\infty\f{s_k}{(-29700)^k}k^3T_k(14,1)T_k(11,-11)^2
\\=&\f{237825864689214465834201682328565\sqrt5}{3584\pi},\end{aligned}
\tag{VIII4S}\]
where
\begin{align*}r_k:=&1921491176698673825280k^2-2877654926621770908536k
\\&+2447787507626240957299
\end{align*}
and
\begin{align*}s_k:=&46821333085900052992621725937920k^2
\\&-71580243731698230860035196766952k
\\&+47098473251985583734378482612735.
\end{align*}
\end{conjecture}

\begin{conjecture} We have
\begin{align*}&\sum_{k=1}^\infty\f{P_1(k)}{3136^k}k^3\bi{2k}kT_k(14,1)T_k(17,16)
\\=&-\f{70662375766828100018060759}{24\pi}\tag{IX1S}
\end{align*}
and
\begin{align*}&\sum_{k=1}^\infty\f{P_2(k)}{3136^k}k^3\bi{2k}kT_k(2,81)T_k(14,81)
\\=&-\f{49}{2916\pi}
(16031198128567258016+6632521068118111355\sqrt5),\tag{IX2S}
\end{align*}
where
\begin{align*}P_1(k):=&8566416619450156417918110k^2-73211973506397665606012003k
\\&+109289868688706912564582106
\end{align*}
and
\begin{align*}P_2(k):=&650373941133830880k^2-11795256486234065124k
\\&+5657835264086004473.
\end{align*}
\end{conjecture}

Those $T_n=T_n(1,1)$ ($n\in\N$) are central trinomial coefficients. Motivated by the author's conjectural identity
$$\sum_{k=1}^\infty\f{(105k-44)T_{k-1}}{k^2\bi{2k}k^23^{k-1}}=6\log3+\f{5\pi}{\sqrt3}$$
(cf. \cite[(10.1)]{S20}), we pose the following conjecture.

\begin{conjecture} We have
$$\sum_{k=1}^\infty\f{(33k^2+32k+8)T_{k-1}(8,-2)}{k(2k+1)^2\bi{2k}k^218^{k-1}}=3\log2$$
and
$$\sum_{k=1}^\infty\f{(15k^3+8k^2-4k-2)T_{k-1}(2,-1)}{k^2(2k+1)^2\bi{2k}k^2}=\f{\pi}{4\sqrt2}.$$
Also,
\begin{gather*}\sum_{k=0}^\infty\f{(32340k^2-28975k-63)T_k}{\bi{2k}k^23^k}=1024+216\log3-\f{388}9\sqrt3\,\pi,
\\\sum_{k=0}^\infty\f{(625k^3-2125k^2+1735k-249)T_k(2,-1)}{\bi{2k}k^2}=\f{9\pi}{\sqrt2}-384,
\\\sum_{k=0}^\infty\f{(11979k^3-50303k^2+43769k-9725)T_k(8,-2)}{\bi{2k}k^218^k}=36(\log2-288).
\end{gather*}
\end{conjecture}

\section{Other series of type S}

The author \cite[(8.1)]{S20} observed the identity
$$\sum_{k=0}^\infty\f{145k+9}{900^k}\beta_kT_k(52,1)=\f{285}{\pi}.$$
Motivated by this we formulate the following conjecture.

\begin{conjecture} {\rm (i)} We have
\begin{align*}\sum_{k=0}^\infty\f{kt_k^{(1)}}
{900^k}\beta_kT_k(52,1)
=&-\f{498507105}{\pi}
\\\sum_{k=0}^\infty\f{k^2t_k^{(2)}}
{900^k}\beta_kT_k(52,1)=&-\f{64232846618769}{\pi},
\\\sum_{k=0}^\infty\f{k^3t_k^{(3)}}
{900^k}\beta_kT_k(52,1)=&-\f{1002142163415074991}{\pi},
\\\sum_{k=0}^\infty\f{(k^2+1)t_k^+}
{900^k}\beta_kT_k(52,1)=&-\f{2512943751517782}{\pi},
\\\sum_{k=0}^\infty\f{(k^2-1)t_k^-}
{900^k}\beta_kT_k(52,1)=&-\f{213755454707472}{\pi},
\end{align*}
where
$$t_k^{(1)}=31642625k^2-170395425k-234814102,$$
\begin{align*}t_k^{(2)}=4936966494550k^2-65077136905125k+134224328497289,
\end{align*}
\begin{align*}t_k^{(3)}=&2822066506655501225k^2-45445823379495456435k
\\&+134916193886322618984,
\end{align*}
\begin{align*}t_k^+=&71236503338400k^2-613367241395375k+111921391484847,
\\t_k^-=&40085552483775k^2-295218039299125k-125724302470098.
\end{align*}

{\rm (ii)} We have
$$\sum_{k=0}^\infty\f{k^j\beta_kT_k(52,1)}{900^k}P_j(k)=0$$
for all $j=0,1,2$, where
\begin{align*}
P_0(k)=&601209875k^3-3237513075k^2+357434077k+299104263,
\\P_1(k)=&3493703960875k^3-48937877104875k^2
\\&+110522556574445k+21410948872923,
\\P_2(k)=&450165200053206075k^3-7261629784552340695k^2
\\&+21683271651287795913k-334047387805024997.
\end{align*}
\end{conjecture}

Motivated by Theorem \ref{Th1.4}, we pose the following two conjectures.

\begin{conjecture} \label{D-Conj} We have
\begin{align*}
\sum_{k=1}^\infty\f{(k-1)(9k+1)}{(-32)^k}k^2\Domb(k)=&\f4{3\pi},
\\\sum_{k=1}^\infty\f{27k^2-12k+17}{(-32)^k}k^3\Domb(k)=&\f2{\pi},
\\\sum_{k=1}^\infty\f{51k^2-60k+89}{(-32)^k}k^4\Domb(k)=&-\f2{3\pi},
\\\sum_{k=1}^\infty\f{801k^2-796k+1371}{(-32)^k}k^5\Domb(k)=&-\f{248}{3\pi},
\\\sum_{k=1}^\infty\f{12339k^2-19308k+46729}{(-32)^k}k^6\Domb(k)=&\f{4886}{3\pi},
\\
\sum_{k=1}^\infty\f{15k^2-51k+22}{64^k}k^3\Domb(k)=&-\f{8\sqrt3}{27\pi},
\\\sum_{k=1}^\infty\f{110k^2-453k+229}{64^k}k^4\Domb(k)=&-\f{8\sqrt3}{45\pi},
\\\sum_{k=1}^\infty\f{1145k^2-7119k+8368}{64^k}k^5\Domb(k)=&\f{7832\sqrt3}{135\pi},
\\\sum_{k=1}^\infty\f{376560k^2-3010921k+4944543}{64^k}k^6\Domb(k)=&\f{3050552\sqrt3}{405\pi}.
\end{align*}
Also,
\begin{align*}
\sum_{k=0}^\infty\f{k\Domb(k)}{(-32)^k}(27k^3-6k^2+9k+2)=&0,
\\\sum_{k=0}^\infty\f{k^2\Domb(k)}{(-32)^k}(54k^3-51k^2+58k+3)=&0,
\\\sum_{k=0}^\infty\f{k^3\Domb(k)}{(-32)^k}(9k^3-9k^2+15k+1)=&0,
\\\sum_{k=0}^\infty\f{k^4\Domb(k)}{(-32)^k}(9k^3-80k^2+99k-124)=&0,
\\\sum_{k=0}^\infty\f{k^5\Domb(k)}{(-32)^k}(1116k^3-319k^2+2808k+2443)=&0,
\\\sum_{k=0}^\infty\f{k^6\Domb(k)}{(-32)^k}(21987k^3-214467k^2+316525k-651453)=&0,
\end{align*}
and
\begin{align*}
\sum_{k=0}^\infty\f{k\Domb(k)}{64^k}(45k^3-75k^2-21k-5)=&0,
\\\sum_{k=0}^\infty\f{k^2\Domb(k)}{64^k}(45k^3-153k^2+71k-3)=&0,
\\\sum_{k=0}^\infty\f{k^3\Domb(k)}{64^k}(25k^3-105k^2+59k-3)=&0,
\\\sum_{k=0}^\infty\f{k^4\Domb(k)}{64^k}(15k^3+377k^2-1827k+979)=&0,
\\\sum_{k=0}^\infty\f{k^5\Domb(k)}{64^k}(132165k^3-1108949k^2+2059839k-381319)=&0,
\\\sum_{k=0}^\infty\f{k^6\Domb(k)}{64^k}(1143957k^3-6418845k^2-10015633k+43834497)=&0.
\end{align*}
\end{conjecture}
\begin{remark} Though we have proved \eqref{D-32}, we are unable to prove the first identity
in Conjecture \ref{D-Conj}. This is because we could not prove the auxilliary identity
\begin{equation}\label{v_k'}\sum_{k=0}^\infty(27k^3-6k^2+9k+2)v_k'=0
\end{equation}
with $v_k'=k\Domb(k)/(-32)^k$. By the Zeilberger algorithm, we find the recurrence
$$(n+1)^4v_n'+n(2n+3)(5n^2+15n+12)v_{n+1}'+16n(n+1)(n+2)^2v_{n+2}'=0$$
from which we have
$$\sum_{k=0}^\infty(27k^4-39k^3+32k^2-3k-1)v_k'=0.$$
Note that this is different from the desired \eqref{v_k'}.
\end{remark}

\begin{conjecture} We have the following identities:
\begin{align*}\sum_{k=1}^\infty\f{1048k^2-5526k+2651}{36^k}k^3f_k^{(4)}=-\f{3618\sqrt{15}}{625\pi},
\tag{Y1S}
\\\sum_{k=1}^\infty\f{304k^2-576k+521}{(-64)^k}k^3f_k^{(4)}=-\f{25696\sqrt{15}}{151875\pi},
\tag{Y2S}
\\\sum_{k=1}^\infty\f{157260600k^2-373237478k+150818747}{196^k}k^3f_k^{(4)}
=-\f{952574\sqrt7}{5\pi},\tag{Y3S}
\end{align*}
\begin{align*}
\sum_{k=1}^\infty\f{1691296k^2-3241287k+1756217}{(-324)^k}k^3f_k^{(4)}=
\f{20434599\sqrt5}{20000\pi},\tag{Y4S}
\end{align*}
\[\sum_{k=1}^\infty\f{833112800k^2 -1793574801k+ 878937679}{1296^k}k^3f_k^{(4)}=-\f{1499005521\sqrt2}{10240\pi},\tag{Y5S}
\]
\[\begin{aligned}&\sum_{k=1}^\infty\f{47808294003072k^2-102482715691400k+52422407372915}{5776^k}k^3f_k^{(4)}
\\&\qquad\quad=-\f{122626206796\sqrt{95}}{625\pi}.
\end{aligned}\tag{Y6S}\]
\end{conjecture}

\begin{conjecture}\label{mixed}
{\rm (i)} We have \begin{align*}\sum_{k=1}^\infty\f{6700k^2-25077k+6239}{96^k}k^3\bi{2k}kf_k=&-\f{5787\sqrt2}{40\pi},
\\\sum_{k=1}^\infty\f{100563606k^2-179847747k+97593215}{(-400)^k}k^3\bi{2k}kf_k
=&\f{642350}{3\pi},
\end{align*}
where $f_k$ denotes the Franel number $\sum_{j=0}^k\bi kj^3$.

{\rm (ii)} We have
\begin{align*}\sum_{k=1}^\infty\f{21k^2-109k-26}{32^k}kT_kZ_k
=&\f{8(250+259\sqrt5)}{375\pi},
\\\sum_{k=1}^\infty\f{1170k^2-14621k+19673}{32^k}k^2T_kZ_k
=&-\f{88(51250+24079\sqrt5)}{375\pi},
\\\sum_{k=1}^\infty\f{885285k^2-14121391k+29366404}{32^k}k^3T_kZ_k
=&\f{8(21456450+3244519\sqrt5)}{125\pi},
\end{align*}
where
$$Z_k=\sum_{j=0}^k\bi kj\bi{2j}j\bi{2(k-j)}{k-j}$$
as introduced by D. Zagier \cite{Zag}.

{\rm (iii)} We have
\begin{align*}\sum_{k=1}^\infty\f{148778208k^2-813461721k+717359335}{(-100)^k}k^3S_k(1,25)
=-\f{19433301932825}{995328\pi},
\end{align*}
where
$$S_k(b,c)=\sum_{j=0}^k\bi kj^2T_j(b,c)T_{k-j}(b,c).$$
\end{conjecture}
\begin{remark} For known $\f1{\pi}$-series involving Franel numbers, see \cite{ChCo,CTYZ}.
Part (ii) of Conjecture \ref{mixed} is motivated by the author's discovery (cf. \cite[(5.1)]{S20})
$$\sum_{k=0}^\infty\f{5k+1}{32^k}T_kZ_k=\f 8{3\pi}(2+\sqrt5),$$
while part (iii) is inspired by the identity
$$\sum_{k=0}^\infty(3k+1)\f{S_k(1,25)}{(-100)^k}=\f{25}{8\pi}$$
obtained by the author (cf. \cite[(1.91)]{S20}).
\end{remark}

\begin{conjecture}\label{324} {\rm (i)} For $n\in\N$ let
$$a_n=\sum_{k=0}^n\bi nk\bi{n+2k}{2k}\bi{2k}k(-324)^{n-k}.$$
Then
\begin{align*}\sum_{k=0}^\infty\f{k\bi{2k}ka_k}{2160^k}(3750642k^2-27417879-200413)=&\f{6028830}{\pi},
\\\sum_{k=0}^\infty\f{k^2\bi{2k}ka_k}{2160^k}(7297838982k^2+44777877957k-249182264449) =&\f{56371645440}{\pi},
\end{align*}
and
\begin{align*}&\sum_{k=0}^\infty\f{k^3\bi{2k}ka_k}{2160^k}(1296246139663698k^2-604571242044753k+ 1636493480946725)
\\&\qquad\quad =-\f{56023236343530}{7\pi}.
\end{align*}
Also,
\begin{align*}\sum_{k=0}^\infty\f{\bi{2k}ka_k}{2160^k}(36414k^3-266193k^2-234124k-66987)=&0,
\\\sum_{k=0}^\infty\f{k\bi{2k}ka_k}{2160^k}(30114378k^3+40060179k^2+29646217k+7732736)=&0,
\end{align*}
and
\begin{align*}
\sum_{k=0}^\infty\f{k^2\bi{2k}ka_k}{2160^k}P(k)=0
\end{align*}
where
$$P(k)=281579848704k^3-131104215018k^2+356871786301k-7684943257.$$

{\rm (ii)} For $n\in\N$ let
$$b_n=\sum_{k=0}^n\bi nk\bi{n+2k}{2k}\bi{2k}k\l(-\f23\r)^{3k}.$$
Then
\begin{align*}\sum_{k=0}^\infty\f{k\bi{2k}kb_k}{(-20)^k}(1440k^2-10980k-559)=&\f3{8\pi}(4255\sqrt6+5904\sqrt{15}),
\\\sum_{k=0}^\infty\f{k^2\bi{2k}kb_k}{(-20)^k}(6192k^2-99894k+37303)=&-\f3{64\pi}(168935\sqrt6+171648\sqrt{15}),
\end{align*}
and
\begin{align*}
&\sum_{k=0}^\infty\f{k^3\bi{2k}kb_k}{(-20)^k}(93108288k^2-1995628400k+ 241469595)
\\&\qquad=\f3{128\pi}(3562267392\sqrt{15}-403810135\sqrt6).
\end{align*}
\end{conjecture}
\begin{remark}  Conjecture \ref{324} is motivated by the author's conjectural series
$$\sum_{k=0}^\infty\f{357k+103}{2160^k}\bi{2k}ka_k=\f{90}{\pi}
\ \t{and}\ \sum_{k=0}^\infty\f{24n+5}{(-20)^k}\bi{2k}kb_k=\f3{2\pi}(5\sqrt6+4\sqrt{15})$$
(cf. \cite[(4.34)]{S14d} and \cite[(5.7)]{S20p}).
\end{remark}

\begin{conjecture}\label{576} For $n\in\N$ let
$$c_n=\sum_{k=0}^n 5^k\f{\bi{2k}k^2\bi{2(n-k)}{n-k}^2}{\bi nk}.$$
Then
\begin{gather*}\sum_{k=0}^\infty\f{k\bi{2k}kc_k}{576^k}(68992k^2-129336k-9691)=\f{27}{2\pi}(2092+1577\sqrt2),
\\\sum_{k=0}^\infty\f{k^2\bi{2k}kc_k}{576^k}(60781952k^2-396958056k+220310929)
\\=-\f{81}{4\pi}(3688232+2390989\sqrt2),
\\\sum_{k=0}^\infty\f{(k^2+1)\bi{2k}kc_k}{576^k}(25180969856k^2 -118523107336k -20012768875)
\\ =-\f{9}{4\pi}(37415939872+18625995101\sqrt2),
\\\sum_{k=0}^\infty\f{k^3\bi{2k}kc_k}{576^k}(197398592384k^2-1384797901272k+855551417353)
\\ =\f{243}{56\pi}(4385109711\sqrt2-1128081200),
\end{gather*}
and
$$\sum_{k=0}^\infty\f{k^4\bi{2k}kc_k}{576^k}Q(k)
=\f{243}{112\pi}(204062865832352+524948955975809\sqrt2),$$
where
$$Q(k)=5366018489638016k^2-50377737067775448k+54028245883253707.$$
\end{conjecture}
\begin{remark} Conjecture \ref{576} is motivated by the author's conjectural series
$$\sum_{k=0}^\infty\f{28k+5}{576^k}\bi{2k}kc_k=\f 9{\pi}(2+\sqrt2)$$
(cf. \cite[(8)]{S14c}).
\end{remark}

As there are many Ramanujan-type series, we cannot list all of their variants of type $S$ here.
We now conclude our paper and hope that our conjectures will stimulate further research.
\medskip

\noindent{\bf Acknowledgment}. The author would like to thank the referee for his/her helpful comments.

\end{document}